\documentclass[10pt,english,french,british]{article}
\usepackage[T1]{fontenc}
\usepackage{array}
\usepackage{amsmath}
\usepackage{graphicx}
\usepackage{amssymb}
\usepackage{esint}
\usepackage[authoryear]{natbib}

\makeatletter

\providecommand{\tabularnewline}{\\}

\usepackage[utf8x]{inputenc}
\usepackage[T1]{fontenc}
\usepackage{longtable}
\usepackage{nicefrac}
\usepackage{bm}
\usepackage{amsfonts}
\usepackage{textcomp}
\usepackage{mathrsfs}
\usepackage{gensymb}
\usepackage{amsthm}
\usepackage{xcolor}
\usepackage{bbm}

 \newtheoremstyle{mestheoremes}{\topsep}{\topsep}{\itshape}{}{\normalfont}{.}{ }{\thmname{#1}\thmnumber{ #2}\thmnote{ (#3)}} 

 \newtheoremstyle{mesdefinitions}{\topsep}{\topsep}{\itshape}{}{\itshape}{.}{ }{\thmname{#1}\thmnumber{ #2}\thmnote{ (#3)}} 

 \newtheoremstyle{mestheoremesp}  {\topsep}{\topsep}{\itshape}   {}{\bfseries}{.}{\newline }{\thmname{#1}\thmnumber{ #2}\thmnote{ : #3}}

\newtheoremstyle{mesdefinitionsp}{\topsep}{\topsep}{\normalfont}{}{\bfseries}{.}{\newline }{\thmname{#1}\thmnumber{ #2}\thmnote{ : #3}}

\theoremstyle{mestheoremes}

\newtheorem{hypo}{A}

\theoremstyle{mestheoremesp}   

\newtheorem{thme}{Theorem}
\newtheorem{prop}[thme]{Proposition}
\newtheorem{lem}[thme]{Lemma}
\newtheorem{defi}[thme]{Definition}

\theoremstyle{mesdefinitionsp}

\theoremstyle{mesdefinitions}
\newtheorem{req}[thme]{Remark}

\newtheorem{result}[thme]{Result}

\newcommand{\norm}[1]{\Vert #1 \Vert}

\newcommand{\rond}[1]{\mathscr{#1}}
\newcommand{\E}{\mathbb{E}}
\newcommand{\Prob}{\mathbb{P}}
\newcommand{\units}[1]{\mathbbmss{1}_{#1}}
\newcommand{\ep}{\varepsilon}

\title{}
\author{Emeline Schmisser\\
Universit\'e Lille 1\\
Laboratoire Paul Painlev\'e\\
\small{emeline.schmisser@math.univ-lille1.fr}}

\makeatother
\usepackage{babel}

\makeatother

\usepackage{babel}
\addto\extrasfrench{\providecommand{\fg}{\ifdim\lastskip>\z@\unskip\fi~\frqq}}

\begin{document}

\title{Non-parametric adaptive estimation of the drift for a jump diffusion
process}
\maketitle
\begin{abstract}
In this article, we consider a jump diffusion process $\left(X_{t}\right)_{t\geq0}$
observed at discrete times $t=0,\Delta,\ldots,n\Delta$. The sampling
interval $\Delta$ tends to 0 and $n\Delta$ tends to infinity. We
assume that $\left(X_{t}\right)_{t\geq0}$ is ergodic, strictly stationary
and exponentially $\beta$-mixing. We use a penalized least-square
approach to compute two adaptive estimators of the drift function
$b$. We provide bounds for the risks of the two estimators. 
\end{abstract}

\section{Introduction}

We consider a general diffusion with jumps: 

\begin{equation}
dX_{t}=b(X_{t})dt+\sigma(X_{t})dW_{t}+\xi(X_{t^{-}})dL_{t}\quad\textrm{and}\quad X_{0}=\eta\label{eq:EDS_jumps}
\end{equation}
where $L_{t}$ is a centred pure jump Levy process: \[
dL_{t}=\int_{z\in\mathbb{R}}z\left(\mu(dt,dz)-dt\nu(dz)\right)\]
with $\mu$ a random Poisson measure with intensity measure $\nu(dz)dt$
such that $\int_{z\in\mathbb{R}}z^{2}\nu(dz)<\infty$. The compensated
Poisson measure $\tilde{\mu}$ is defined by $\tilde{\mu}(dt,dz)=\mu(dt,dz)-\nu(dz)dt.$
The random variable $\eta$ is independent of $(W_{t},L_{t})_{t\geq0}$.
Moreover, $(W_{t})_{t\geq0}$ and $(L_{t})_{t\geq0}$ are independent. 

This process is observed with high frequency (at times $t=0,\Delta,\ldots,n\Delta$
where, as $n$ tends to infinity, the sampling interval $\Delta\rightarrow0$ and the time of observation $n\Delta\rightarrow\infty$). It is
assumed to be ergodic, stationary and exponentially $\beta$-mixing
(see \citet{masuda2007} for sufficient conditions). Our aim is to
construct a non-parametric estimator of  $b$ on a compact set $A$.

The non-parametric estimation of $b$ and $\sigma$ for a diffusion
process observed with high-frequency is well-known (see for instance
\citet{hofmann99} and \citet{comtegenon2007}). Diffusion processes
with jumps are used in various fields, for instance in finance, for
modelling the growth of a population, in hydrology, in medical science,
$\ldots$, but there exist few results for the non-parametric estimation
of $b$ and $\sigma$. \citet{mai2012} and  \citet{shimizu_yoshida2006} construct maximum-likelihood
estimators of parameters of $b$. Their estimators reach
the standard rate of convergence, $\sqrt{n\Delta}$. \citet{shimizu2008}
and \citet{mancinireno2011} use a kernel estimator to obtain non
parametric threshold estimators of $\sigma$. \citet{mancinireno2011}
also construct a non-parametric truncated estimator of $b$, but
only when $L_{t}$ is a compound Poisson process. To our knowledge,
minimax rates of convergences for non-parametric estimators of $b$,
$\sigma$ or $\xi$ for jump-diffusions processes are not available in the literature 
(see \cite{hofmann99} or \cite{gobethoffmannreiss2004} for rates of convergence 
for diffusions processes). 

In this paper, we use model selection to construct two non-parametric
estimators of $b$ under the asymptotic framework $\Delta\rightarrow0$
and $n\Delta\rightarrow\infty$. This method was introduced by\foreignlanguage{french}{
\citet{birgemassart98}}. 

First, we introduce a sequence of linear subspaces $S_{m}\subseteq L^{2}(A)$
and, for each $m$, we construct an estimator $\hat{b}_{m}$ of $b$
by minimising on $S_{m}$ the contrast function: \[
\gamma_{n}(t)=\frac{1}{n}\sum_{k=1}^{n}\left(Y_{k\Delta}-t(X_{k\Delta})\right)^{2}\quad\textrm{where}\quad Y_{k\Delta}=\frac{X_{(k+1)\Delta}-X_{k\Delta}}{\Delta}.\]
We obtain a collection of estimators of the drift function $b$ and
we bound their risks (Theorem \ref{thme_risque_m_fixe}). Then, we
introduce a penalty function to select the {}``best'' dimension
$m$ and we deduce an adaptive estimator $\hat{b}_{\hat{m}}$. Under
the assumption that $\nu$ is sub-exponential, that is if there exist
two positive constants $C$, $\lambda$ such that, for $z$ large
enough, $\nu([-z,z]^{c})\leq Ce^{-\lambda z}$, the risk bound of
$\hat{b}_{\hat{m}}$ is exactly the same as for a diffusion without
jumps (Theorem \ref{thme_adaptatif}) (see \citet{comtegenon2007}
or \citet{hofmann99}).

In a second part, we do not assume that $\nu$ is sub-exponential and
we construct a truncated estimator $\tilde{b}_{m}$ of $b$. We minimise
the contrast function \[
\tilde{\gamma}_{n}(t)=\frac{1}{n}\sum_{k=1}^{n}\left(Y_{k\Delta}\units{\left|Y_{k\Delta}\right|\leq C_{\Delta}}-t(X_{k\Delta})\right)^{2}\quad\textrm{where}\quad C_{\Delta}\propto\sqrt{\Delta}\ln(n)\]
in order to obtain a new estimator $\tilde{b}_{m}$. As in the first
part, we introduce a penalty function to obtain an adaptive estimator
$\tilde{b}_{\tilde{m}}$. The risk bound of this adaptive estimator
depends on the Blumenthal-Getoor index of $\nu$ (Theorems \ref{thme_sauts_coupes_m_fixe }
and \ref{thme_sauts_coupes_adaptif}). 

In Section \ref{sec:Model}, we present the model and its assumptions.
In Sections \ref{sec:Estimation} and \ref{sub:Cutted-jumps}, we
construct the estimators and bound their risks. Some simulations are
presented in Section \ref{sec:simulation}. Proofs are gathered in
Section \ref{sec:Proofs}.

\section{Assumptions \label{sec:Model}}

\subsection{Assumptions on the model }

We consider the following assumptions:

\begin{hypo}\label{hypo_eds}

The functions $b$, $\sigma$ and $\xi$ are Lipschitz.

\end{hypo}

\begin{hypo}\label{hypo_densite}
\begin{enumerate}
\item The function $\sigma$ is bounded from below and above: \[
\exists\sigma_{0},\sigma_{1},\;\forall x\in\mathbb{R},\quad0<\sigma_{1}\leq\sigma(x)\leq\sigma_{0}.\]

\item The function $\xi$ is bounded: $\exists\xi_{0},\;\forall x\in\mathbb{R},\quad0\leq\xi(x)\leq\xi_{0}.$
\item The drift function $b$ is elastic: there exists a constant $M$ such
that, for any $x\in\mathbb{R}$, $\left|x\right|>M$: $xb(x)\lesssim-\left|x\right|^{2}.$
\item The Lévy measure $\nu$ satisfies: \[
\nu(\{0\})=0,\quad\int_{-\infty}^{\infty}z^{2}\nu(dz)=1\quad\textrm{and}\quad\int_{-\infty}^{\infty}z^{4}\nu(dz)<\infty.\]

\end{enumerate}
\end{hypo}

Under Assumption A\ref{hypo_eds}, the stochastic differential equation
\eqref{eq:EDS_jumps} admits a unique strong solution. According to
\citet{masuda2007}, under Assumptions A\ref{hypo_eds} and A\ref{hypo_densite},
the process $\left(X_{t}\right)$ admits a unique invariant probability
$\varpi$ and satisfies the ergodic theorem: for any measurable function
$g$ such that $\int\vert g(x)\vert\varpi(dx)<\infty$, when $T\rightarrow\infty$,
\[
\frac{1}{T}\int_{0}^{T}g(X_{s})ds\rightarrow\int g(x)\varpi(dx).\]
This distribution has moments of order 4. Moreover, \citet{masuda2007}
also ensures that under these assumptions, the process $\left(X_{t}\right)$
is exponentially $\beta$-mixing. Furthermore, if there exist two constants
$c$ and $n_{0}$ such that, for any $x\in\mathbb{R}$, $\xi^{2}(x)\geq c(1+\vert x\vert)^{-n_{0}}$,
then \citet{ishikawakunita2006} ensure that a smooth transition density
exists. 

\begin{hypo}\label{hypo_stationnarite}
\begin{enumerate}
\item The stationary measure $\varpi$ admits a density $\pi$ which is
bounded from below and above on the compact interval $A$: \[
\exists\pi_{0},\pi_{1},\;\forall x\in A,\quad0<\pi_{1}\leq\pi(x)\leq\pi_{0}.\]

\item The process $\left(X_{t}\right)_{t\geq0}$ is stationary ($\eta\sim\varpi(dx)=\pi(x)dx$).
\end{enumerate}
\end{hypo}
The first part of this assumption is automatically satisfied if $\xi=0$ (that is if $(X_t)_{t\geq 0}$ is a diffusion process). 
The following proposition is very useful for the proofs. It is derived from Result \ref{res_burkholder}. 

\begin{prop}\label{prop_gloter}

Under Assumptions A\ref{hypo_eds}-A\ref{hypo_stationnarite}, for
any $p\geq1$, there exists a constant $c(p)$ such that\emph{, }if
$\int_{\mathbb{R}}z^{2p}\nu(dz)<\infty$: \[
\E\left(\sup_{s\in[t,t+h]}\left(X_{s}-X_{t}\right)^{2p}\right)\leq c(p)h.\]

\end{prop}

\subsection{Assumptions on the approximation spaces}

In order to construct an adaptive estimator of $b$, we use model
selection: we compute a collection of estimators $\hat{b}_{m}$ of
$b$ by minimising a contrast function $\gamma_{n}(t)$ on a vectorial
subspace $S_{m}\subset L^{2}(A)$, then we choose the best possible
estimator using a penalty function $pen(m)$. The collection of vectorial
subspaces $\left(S_{m}\right)_{m\in\rond{M}_{n}}$ has to satisfy
the following assumption: 

\begin{hypo}\label{hypo_espaces}$\ $
\begin{enumerate}
\item The subspaces $S_{m}$ have finite dimension $D_{m}$.
\item The sequence of vectorial subspaces $(S_{m})_{m\geq0}$ is increasing:
for any $m,$ $S_{m}\subseteq S_{m+1}$. 
\item Norm connexion: there exists a constant $\phi_{1}$ such that, for
any $m\geq0$, any $t\in S_{m}$, \[
\left\Vert t\right\Vert _{\infty}^{2}\leq\phi_{1}D_{m}\left\Vert t\right\Vert _{L^{2}}^{2}\]
where $\norm{.}_{L^{2}}$ is the $L^{2}$-norm and $\norm{.}_{\infty}$
is the sup-norm on $A$. 
\item For any $m\in\mathbb{N}$, there exists an orthonormal basis $\left(\psi_{\lambda}\right)_{\lambda\in\Lambda_{m}}$
of $S_{m}$ such that \[
\forall\lambda,\quad\textrm{\textrm{card}}\left(\lambda',\,\left\Vert \psi_{\lambda}\psi_{\lambda'}\right\Vert _{\infty}\neq0\right)\leq\phi_{0}\]
where $\phi_{0}$ does not depend on $m$. 
\item For any function $t$ belonging to the unit ball of the Besov space
$\rond{B}_{2,\infty}^{\alpha}$, \[
\exists C,\:\forall m\quad\left\Vert t-t_{m}\right\Vert _{L^{2}}^{2}\leq CD_m^{-2\alpha}\]
where $t_{m}$ is the $L^{2}$ orthogonal projection of $t$ on $S_{m}$. 
\end{enumerate}
\end{hypo}

The subspaces generated by piecewise polynomials, compactly supported
wavelets or spline functions satisfy A\ref{hypo_espaces} (see \citet{devorelorentz}
and \citet{meyer} for instance).

\section{Estimation of the drift \label{sec:Estimation}}

By analogy with \citet{comtegenon2007}, we decompose $Y_{k\Delta}$
in the following way: 
\begin{equation}
Y_{k\Delta}=\frac{X_{(k+1)\Delta}-X_{k\Delta}}{\Delta}=b(X_{k\Delta})+I_{k\Delta}+Z_{k\Delta}+T_{k\Delta}
\label{eq:def_Y}
\end{equation}
where 
\begin{gather*}
I_{k\Delta}=\frac{1}{\Delta}\int_{k\Delta}^{(k+1)\Delta}\left(b(X_{s})-b(X_{k\Delta})\right)ds,
\quad 
Z_{k\Delta}=\frac{1}{\Delta}\int_{k\Delta}^{(k+1)\Delta}\sigma(X_{s})dW_{s}\\
T_{k\Delta}=\frac{1}{\Delta}\int_{k\Delta}^{(k+1)\Delta}\xi(X_{s^{-}})dL_{s}.
\end{gather*}
The terms $Z_{k\Delta}$ and $T_{k\Delta}$ are martingale increments.
Let us introduce the mean square contrast function 
\begin{equation}
\gamma_{n}(t)=\frac{1}{n}\sum_{k=1}^{n}\left(Y_{k\Delta}-t\left(X_{k\Delta}\right)\right)^{2}\label{eq:def_gamma}.
\end{equation}
We can always minimise $\gamma_n(t)$ on $S_m$, but the minimiser may be not unique. That is why we introduce the 
 empirical risk
 \begin{equation}
\rond{R}_{n}(t)=\E\left(\left\Vert t-b_{A}\right\Vert _{n}^{2}\right)
\quad\textrm{where}\quad
\left\Vert t\right\Vert _{n}^{2}=\frac{1}{n}\sum_{k=1}^{n}t^{2}\left(X_{k\Delta}\right)
\quad\textrm{and}\quad 
t_{A}=t\units{A}.
\label{eq:def_risque}
\end{equation}
We consider the asymptotic framework: 
\[
\Delta\rightarrow0,\quad n\Delta\rightarrow\infty.
\]
For any 
$m\in\rond{M}_{n}=\left\{ m,\; D_{m}\leq\rond{D}_{n}\right\} $
where 
$\rond{D}_{n}^{2}\leq n\Delta/\ln^{2}(n)$, we construct the
regression-type estimator: 
\[
\hat{b}_{m}=\arg\min_{t\in S_{m}}\gamma_{n}(t).
\]

\begin{thme}\label{thme_risque_m_fixe}

Under Assumptions \emph{A\ref{hypo_eds}-A\ref{hypo_espaces},} the
risk of the estimator with fixed $m$ satisfies: 

\[
\rond{R}_{n}(\hat{b}_{m})\leq3\pi_{1}\left\Vert b_{m}-b_{A}\right\Vert _{L^{2}}^{2}+48(\sigma_{0}^{2}+\xi_{0}^{2})\frac{D_{m}}{n\Delta}+c\Delta
\]
where $b_{m}$ is the orthogonal ($L^{2}$) projection of $b_{A}$
over the vectorial subspace $S_{m}$. The constant $c$ is independent
of $m$, $n$ and $\Delta$. 

\end{thme}

Except for the constant $(\sigma_{0}^{2}+\xi_{0}^{2})$ in the variance
term, this is exactly the bound of the risk that \citet{comtegenon2007}
found for a diffusion process without jumps. 

The bias term, $\left\Vert b_{m}-b_{A}\right\Vert _{L^{2}}^{2}$,
decreases when the dimension $D_{m}$ increases whereas the variance
term $(\sigma_{0}^{2}+\xi_{0}^{2})D_{m}/(n\Delta)$ is proportional
to the dimension. Under the classical assumption $n\Delta^{2}=O(1)$,
the remainder term $\Delta$ is negligible. Thus we need to find a
good compromise between the bias and the variance term. 

\begin{req}

If the regularity of the drift function is known, that is, if $b$
belongs to a ball of a Besov space $\rond{B}_{2,\infty}^{\alpha}$,
then the bias term $\left\Vert b_{m}-b_{A}\right\Vert _{L^{2}}^{2}$
is smaller than $D_{m}^{-2\alpha}$. The best estimator is obtained when the bias term, $\norm{b_m-b_A}_{L^2}^{2}$, and the variance term, $cD_m(n\Delta)^{-1}$, 
are equal, that is 
for $D_{m_{opt}}=\left(n\Delta\right)^{1/(1+2\alpha)}$. In that case, the estimator
risk satisfies: \[
\rond{R}_{n}(\hat{b}_{m_{opt}})\lesssim\left(n\Delta\right)^{-2\alpha/(2\alpha+1)}+\Delta.\]

\end{req}

Let us introduce a penalty function $pen$ such that : \[
pen(m)\geq\kappa(\sigma_{0}^{2}+\xi_{0}^{2})\frac{D_{m}}{n\Delta}\]
and set: \[
\hat{m}=\arg\min_{m\in\rond{M}_{n}}\left\{ \gamma_{n}(\hat{b}_{m})+pen(m)\right\} .\]
We will chose $\kappa$ later. We denote by $\hat{b}_{\hat{m}}$ the
resulting estimator. To bound the risk of the adaptive estimator,
an additional assumption is needed:

\begin{hypo}\label{hypo_loi_xi}
\begin{enumerate}
\item The Lévy measure $\nu$ is symmetric or the function $\xi$  is constant. 
\item The Lévy measure $\nu$ is sub exponential: there exist $\lambda,C>0$
such that, for any $\left|z\right|>1$, $\nu(]-z,z[^{c})\leq Ce^{-\lambda\left|z\right|}$. 
\end{enumerate}
\end{hypo}

\begin{thme}\label{thme_adaptatif}

Under Assumptions A\ref{hypo_eds}-A\ref{hypo_loi_xi}, there exists
a constant $\kappa$ (depending only on $\nu$) such that, if $\rond{D}_{n}^{2}\leq n\Delta/\ln^{2}(n)$:
\[
\mathbb{E}\left(\left\Vert \hat{b}_{\hat{m}}-b_{A}\right\Vert _{n}^{2}\right)\lesssim\inf_{m\in\rond{M}_{n}}\left(\left\Vert b_{m}-b_{A}\right\Vert _{L^{2}}^{2}+pen(m)\right)+\left(\Delta+\frac{1}{n\Delta}\right).\]

\end{thme}

\begin{req}
We can bound $\kappa$ theoretically, however, this bound is in practice
too large for the simulations. In Section 5, we calibrate $\kappa$
by simulations (see \citet{comtegenon2007} for instance). 
If $\sigma$ and $\xi$ are unknown, it is possible to replace them by rough estimators (in fact, we only need upper bounds of $\sigma_0^2$ and $\xi_0^2$). 
It is also possible to performe a completely data-driven calibration of the parameters of the penalty (see \citet{arlotmassart2009}). 
\end{req}

\section{Truncated estimator of the drift  \label{sub:Cutted-jumps}}

Truncated estimators are widely used for the estimation of the diffusion
coefficient of a jump diffusion (see for instance \citet{mancinireno2011}, \citet{shimizu2008} and \citet{mai2012}). Our aim is to construct an adaptive estimator
of $b$ even if Assumption A\ref{hypo_loi_xi} is not fulfilled. To
this end, we cut off the big jumps. Let us introduce the set \[
\Omega_{X,k}=\left\{ \omega,\;\left|X_{(k+1)\Delta}-X_{k\Delta}\right|\leq C_{\Delta}\right\} \]
\foreignlanguage{english}{where $C_{\Delta}=(b_{max}+3)\Delta+\left(\sigma_{0}+4\xi_{0}\right)\sqrt{\Delta}\ln(n)$}
(with $b_{max}=\sup_{x\in A}\left|b(x)\right|$). Let us consider
the random variables \[
\tilde{Y}_{k\Delta}=\frac{X_{(k+1)\Delta}-X_{k\Delta}}{\Delta}\units{\Omega_{X,k}}\units{X_{k\Delta}\in A}.\]
We recall here the definition of the Blumenthal-Getoor index: 

\begin{defi}

The Blumenthal-Getoor index of a Lévy measure is \[
\beta=\inf\left\{ \alpha\geq0,\;\int_{\vert z\vert\leq1}\vert z\vert^{\alpha}\nu(dz)<\infty\right\} .\]

\end{defi}

A compound Poisson process has $\beta=0$.

We assume that the following assumption is fulfilled. 

\begin{hypo}\label{hypo_nu_reguliere}
\begin{enumerate}

\item For $\left|x\right|$ small, $\nu(dx)$ is absolutely continuous with
respect to the Lebesgue measure ($\nu(dx)=n(x)dx$) and: 
\[
\exists\beta\in[0,2[,\exists a_{0},\forall x\in[-a_{0},a_{0}],\quad n(x)\leq Cx^{-\beta-1}.
\]
This implies that the Blumenthal-Getoor index is equal to $\beta$. 
\item The Lévy measure $\nu(z)$ is symmetric for $z$ small: 
\[ 
\exists a_1<a_0, \forall z\in[-a_1,a_1], n(z)=n(-z)
 \]
\item The function $\xi$ is bounded from below: there exists $\xi_{1}>0$
such that, for any $z\in\mathbb{R}$, $0<\xi_{1}\leq\xi(z)$. 
\item The functions $\sigma$ and $\xi$ are $\rond{C}^{2}$, $\xi'$
and $\sigma'$ are Lipschitz. 
\end{enumerate}
\end{hypo}

We consider the following asymptotic framework: 
\[
\frac{n\Delta}{\ln^{2}(n)}\rightarrow \infty,\quad\Delta^{1-\beta/2}\ln^{2}(n)\rightarrow0.
\]
The truncated estimator $\tilde{b}_{m}$ is obtained by minimising
the contrast function: 
\[
\tilde{b}_{m}=\arg\min_{t\in S_{m}}\tilde{\gamma}_{n}(t)
\quad\textrm{where}\quad
\tilde{\gamma}_{n}(t)=\frac{1}{n}\sum_{k=1}^{n}\left(\tilde{Y}_{k\Delta}-t(X_{k\Delta})\right)^{2}.
\]

\begin{thme}[Risk of the non adaptive truncated estimator]\label{thme_sauts_coupes_m_fixe }

Under Assumptions A\ref{hypo_eds}-A\ref{hypo_espaces} and A\ref{hypo_nu_reguliere},
for any $m$ such that $D_{m}\leq\rond{D}_{n}$ where $\rond{D}_{n}^{2}\leq n\Delta/\ln^{2}(n)$: 

\[
\E\left(\left\Vert \tilde{b}_{m}-b_{A}\right\Vert _{n}^{2}\right)\lesssim\left\Vert b_{m}-b_{A}\right\Vert _{L^{2}}^{2}+(\sigma_{0}^{2}+c\Delta^{1/2-\beta/4})\frac{D_{m}}{n\Delta}+\Delta^{1-\beta/2}\ln^{2}(n)+\frac{1}{n\Delta}.\]

\end{thme}
The variance term  is smaller than for the first estimator, but the remainder term 
 depends on the Blumenthal-Getoor index and is
larger than for the first estimator. This remainder term is due to the fact that $\tilde{Y}_{k\Delta}=0$ every time
 $\vert X_{(k+1)\Delta}-X_{k\Delta}\vert>C_{\Delta}$: then
\[\left\vert\E\left(\tilde{Y}_{k\Delta}-b(X_{k\Delta})\right)\right\vert>\left\vert\E\left(Y_{k\Delta}-b(X_{k\Delta})\right)\right\vert.\]
 If $L_{t}$ is
a compound Poisson process, (which implies $\beta=0$) or if $\Delta$ is small enough (see Remark \ref{req_haute_frequence_comparaison_estimateurs}),
 we obtain  a better  inequality than for the non-truncated estimator.

\begin{req}

If $\nu$ is not absolutely continuous, we can prove the weaker inequality:
\[
\E\left(\left\Vert \tilde{b}_{m}-b_{A}\right\Vert _{n}^{2}\right)\lesssim\left\Vert b_{m}-b_{A}\right\Vert _{L^{2}}^{2}+(\sigma_{0}^{2}+\xi_{0}^{2})\frac{D_{m}}{n\Delta}+\Delta^{1-\beta}\ln^{2}(n)+\frac{1}{n\Delta}.\]
In that case, $\tilde{b}_{m}$ converges towards $b_{A}$ only if
$\beta<1$, which implies that $\nu$ has finite variation ($\int_{\mathbb{R}}\vert z\vert\nu(dz)<\infty)$. 
See Remark \ref{ref:nu_pas_continue}. 
\end{req}

\begin{req}\label{req_haute_frequence_comparaison_estimateurs}
Assume that $b_{A}$ belongs to the Besov space $\rond{B}_{2,\infty}^{\alpha}$
and that ${\left\Vert b_A\right\Vert _{\rond{B}_{2,\infty}^{\alpha}}\leq1}$.
The bias-variance compromise $\left\Vert b_{m}-b_{A}\right\Vert _{L^{2}}^{2}+D_{m}/n\Delta$
is minimal when ${m=\log_{2}(n\Delta)/(1+2\alpha)}$, and the risk
satisfies: \[
\E\left(\left\Vert \tilde{b}_{m}-b_{A}\right\Vert _{n}^{2}\right)\lesssim\left(n\Delta\right)^{-2\alpha/(1+2\alpha)}+\Delta^{1-\beta/2}\ln^{2}(n)\]

Let us set $\Delta\sim n^{-\gamma}$ with $\gamma>0$. We have the
following convergence rates:

\begin{center}
\begin{tabular}{|c|c|c|}
\hline 
$\gamma$ & first estimator & truncated estimator\tabularnewline
\hline 
$0<\gamma\leq\frac{2\alpha}{4\alpha+1}\leq\frac{1}{2}$ & $\Delta$ & $\Delta^{1-\beta/2}\ln^{2}(n)$\tabularnewline
\hline 
$\frac{2\alpha}{4\alpha+1}\leq\gamma\leq\frac{2\alpha}{4\alpha+1-\beta\alpha-\beta/2}\leq\frac{1}{2(1-\beta/4)}$ & $\left(n\Delta\right)^{-2\alpha/(2\alpha+1)}$ & $\Delta^{1-\beta/2}\ln^{2}(n)$\tabularnewline
\hline 
$\frac{2\alpha+1}{4\alpha+1-\beta\alpha-\beta/2}\leq\gamma<1$ & $\left(n\Delta\right)^{-2\alpha/(2\alpha+1)}$ & $\left(n\Delta\right)^{-2\alpha/(2\alpha+1)}$\tabularnewline
\hline
\end{tabular}
\par\end{center}

If we have sufficiently high frequency data ($n\Delta^{2(1-\beta/4)}=O(1)$),
then the rate of convergence is $(n\Delta)^{2\alpha/(2\alpha+1}$ for the two estimators. The estimator of \citet{mai2012} converges with the corresponding 
parametric rate, $n\Delta$, if $n\Delta^{3/2-\gamma}=o(1)$ for $\gamma\in]0,1/2[$. 

\end{req}

To construct the adaptive estimator, we use the same penalty function
as in the previous section: \[
pen(m)\geq\kappa\left(\sigma_{0}^{2}+\xi_{0}^{2}\right)\frac{D_{m}}{n\Delta}\]
 and define the adaptive estimator: \[
\tilde{m}=\arg\min_{m\in\rond{M}_{n}}\left\{ \tilde{\gamma}_{n}(\tilde{b}_{m})+pen(m)\right\} .\]

\begin{thme}[Risk of the  adaptive truncated estimator]\label{thme_sauts_coupes_adaptif}

If Assumptions A\ref{hypo_eds}-A\ref{hypo_espaces} and A\ref{hypo_nu_reguliere}
are satisfied, then there exists $\kappa$ such that, if $\rond{D}_{n}^{2}\leq n\Delta/\ln^{2}(n)$: 

\[
\E\left(\left\Vert \tilde{b}_{\tilde{m}}-b_{A}\right\Vert _{n}^{2}\right)\lesssim\min_{m\in\rond{M}_{n}}\left(\left\Vert b_{m}-b_{A}\right\Vert _{n}^{2}+pen(m)\right)+\Delta^{1-\beta/2}\ln^{2}(n)+\frac{1}{n\Delta}.\]

\end{thme}

The adaptive estimator $\tilde{b}_{\tilde{m}}$ automatically realises
the bias/variance compromise.

\section{Numerical simulations and examples \label{sec:simulation}}

\subsection{Models}

We consider the stochastic differential equation: \[
dX_{t}=b(X_{t})dt+\sigma(X_{t})dW_{t}+\xi(X_{t^{-}})dL_{t}\]
 where $L_{t}$ is a compound Poisson process of intensity 1: $L_{t}=\sum_{j=1}^{N_{t}}\zeta_{i}$,
with $N_{t}$ a Poisson process of intensity 1 and $(\zeta_{1},\ldots,\zeta_{n})$
are independent and identically distributed random variables independent
of $(N_{t})$. We denote by $f$ the probability law of $\zeta_{i}$
.

\subsubsection*{Model 1: }

\[
b(x)=-2x,\quad\sigma(x)=\xi(x)=1\quad\textrm{and}\quad f(dz)=\nu(dz)=\frac{1}{2}\delta_{1}+\frac{1}{2}\delta_{-1}.\]

\subsubsection*{Model 2: }

\[
b(x)=-\left(x-1/4\right)^{3}-\left(x+1/4\right)^{3},\quad\sigma(x)=\xi(x)=1\quad\textrm{and}\quad f(dz)=\nu(dz)=\frac{e^{-\lambda\left|z\right|}dz}{2}.\]
We can remark that the function $b$ is not Lipschitz and therefore
does not satisfy Assumption A\ref{hypo_eds}.

\subsubsection*{Model 3: }

We consider the stochastic process of parameters \[
b(x)=-2x+\sin(3x),\quad\sigma(x)=\xi(x)=\sqrt{\frac{3+x^{2}}{1+x^{2}}}\]
and \[
f(dz)=\nu(dz)=\frac{1}{4}\sqrt{\frac{\sqrt{24}}{\left|z\right|}}e^{-\sqrt{\sqrt{24}\left|z\right|}}dz.\]
 Let us remark that $\nu=f$ is not sub-exponential and does not satisfy
A\ref{hypo_loi_xi}. Nevertheless, this model satisfies all the assumptions of Theorem \ref{thme_sauts_coupes_adaptif}.

\subsubsection*{Model 4:}
In this model, the Lévy process is not a compound Poisson process. We set
\[\nu(dz)=\sum_{k=0}^{\infty} 2^{k+2}(\delta_{1/2^{k}}+\delta_{-1/2^{k}}), \quad b(x)=-2x\quad \textrm{and}\quad \sigma(x)=\xi(x)=1.\]
The Blumenthal-Getoor index of this process is such that $\beta>1$.

\subsection{Simulation algorithm (Compound Poisson case)}

We estimate $b$ on the compact interval $A=[-1,1]$.
\begin{enumerate}
\item Simulate random variables $\left(X_{0},X_{\Delta},\ldots,X_{n\Delta}\right)$
thanks to a Euler scheme with sampling interval $\delta=\Delta/5$.
To this end, we use the same simulation scheme as \citet{Rubenthaler_HDR}.
We simulate the times of the jumps ($\tau_{1},\ldots,\tau_{N},\tau_{N+1})$
with $\tau_{N}<n\Delta\leq\tau_{N+1}$ and we fix $X_{0}=0$. \\
If $\delta<\tau_{1}$, we compute 
\[
X_{\delta}=\delta b(X_{0})+\sqrt{\delta}\sigma(X_{0})N
\quad\textrm{with}\quad
 N\sim\rond{N}(0,1).
\]
If $\tau_{1}<\delta$, we first compute 
\[
X_{\tau_{1}}=\tau_{1}b(X_{0})+\sqrt{\tau_{1}}\sigma(X_{0})N+\xi(X_{0})\zeta_{1}
\]
with $N\sim\rond{N}(0,1)$ and $\zeta_{1}\sim f$ is independent of
$N$. If $\delta<\tau_{2}$, we compute 
\[
X_{\delta}=(\delta-\tau_{1})b(X_{\tau_{1}})+\sqrt{\delta-\tau_{1}}\sigma(X_{\tau_{1}})N'
\]
 else we compute 
\[
X_{\tau_{2}}=(\tau_{2}-\tau_{1})b(X_{\tau_{1}})+\sqrt{\tau_{2}-\tau_{1}}\sigma(X_{\tau_{1}})N'+\xi(X_{\tau_{1}})\zeta_{2}
\]
 where $N'\sim\rond{N}(0,1)$ and $\zeta_{2}$ has distribution $f$.
$N$, $N'$, $\zeta_{1}$ and $\zeta_{2}$ are independent. \\

\item Construct the random variables 
\[
Y_{k\Delta}=\frac{X_{(k+1)\Delta}-X_{k\Delta}}{\Delta}
\quad\textrm{and}\quad
\tilde{Y}_{k\Delta}=\frac{X_{(k+1)\Delta}-X_{k\Delta}}{\Delta}\units{\Omega_{X,k}}\units{X_{k\Delta}\in A}.
\]

\item We consider the vectorial subspaces $S_{m,r}$ generated by the spline
functions of degree $r$ (see for instance \citet{schmisser09a}).
In that case $D_{m,r}=\dim(S_{m,r})=2^{m}+r$. For $r\in\left\{ 1,2,3\right\} $
and $m\in\rond{M}_{n}(r)=\{m,D_{m,r}\leq\rond{D}_{n}\}$, we compute
the estimators $\hat{b}_{m,r}$ and $\tilde{b}_{m,r}$ by minimising
the contrast functions $\gamma_{n}$ and $\tilde{\gamma}_{n}$ on
the vectorial subspaces $S_{m,r}$. 
\item For the estimation algorithm, we make a selection of $m$ and $r$
as follows. Using the penalty function 
$pen(m,r):=pen(m)=\kappa(\sigma_{0}^{2}+\xi_{0}^{2})(2^{m}+r)/n\Delta$,
we select the adaptive estimators 
$\hat{b}_{\hat{m},r}$ and $\tilde{b}_{\tilde{m},r}$,
and then choose the best $r$ by minimizing 
$\gamma_{n}(\hat{b}_{\hat{m},r})+pen(\hat{m},r)$
and 
$\tilde{\gamma}_{n}(\tilde{b}_{\tilde{m},r})+pen(\tilde{m},r).$
\end{enumerate}

To calibrate $\kappa$, we run a various number of simulations for a model with known parameters and 
let $\kappa$ vary. When $\kappa$ is too small, the value of m selected by the estimation procedure is in general very high (often maximal). 
When $\kappa$ is too big, the estimator is always linear even if the true function is not. 
We used the true value of $\sigma_{0}^{2}$ and $\xi_{0}^{2}$.

\subsection{Results}

\begin{figure}[p]
\caption{Model 1: Ornstein-Uhlenbeck and binomial law}

\[
b(x)=-2x,\;\sigma(x)=\xi(x)=1\:\textrm{and binomial law}
\]

\label{Flo:figure_model 1}

\begin{centering}
\includegraphics[scale=0.5]{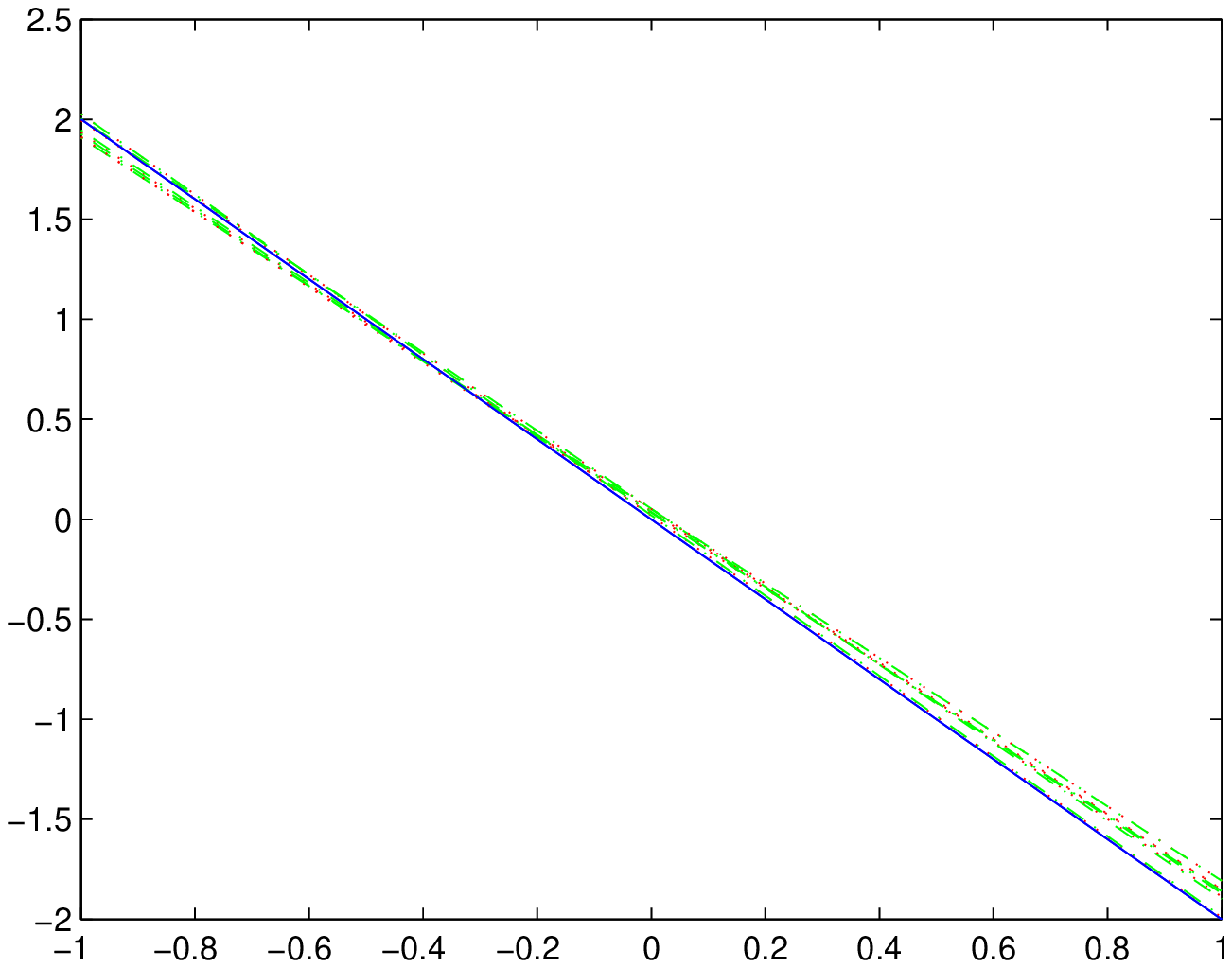}
\par\end{centering}

\centering{}-- : true function $\quad$-.-: first estimator $\quad\ldots$:
truncated estimator \\
$n=10^{4}$ et $\Delta=10^{-1}$
\end{figure}

\begin{figure}[p]
\caption{Model 2: Double well and Laplace law}

\[
b(x)=-\left(x-1/4\right)^{3}-\left(x+1/4\right)^{3},\quad\sigma=\xi=1\textrm{ and Laplace law}
\]

\label{Flo:figure_Model 2}

\centering{}\includegraphics[scale=0.5]{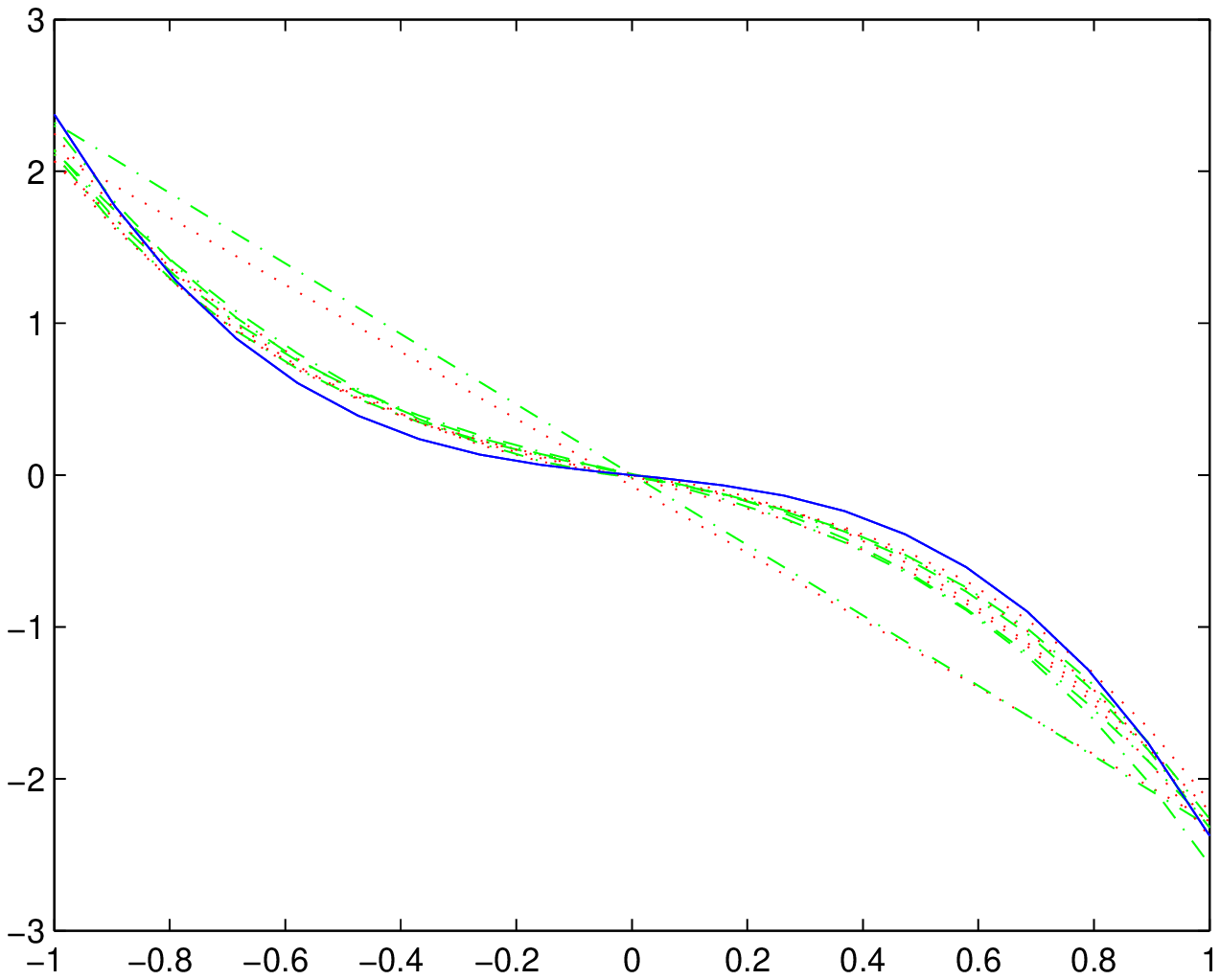}\\
-- : true function$\quad$-.-: first estimator$\qquad\ldots$:
truncated estimator \\
$n=10^{4}$ et $\Delta=10^{-1}$
\end{figure}

\begin{figure}[p]
\caption{Model 3: Sine function }

\begin{centering}
\[
b(x)=-2x+\sin(3x),\;\sigma(x)=\xi(x)=\sqrt{(3+x^{2})/(1+x^{2})}\textrm{ jumps not sub-exponential }
\]
\label{Flo:figure_Model 3}
\par\end{centering}

\begin{centering}
\includegraphics[scale=0.5]{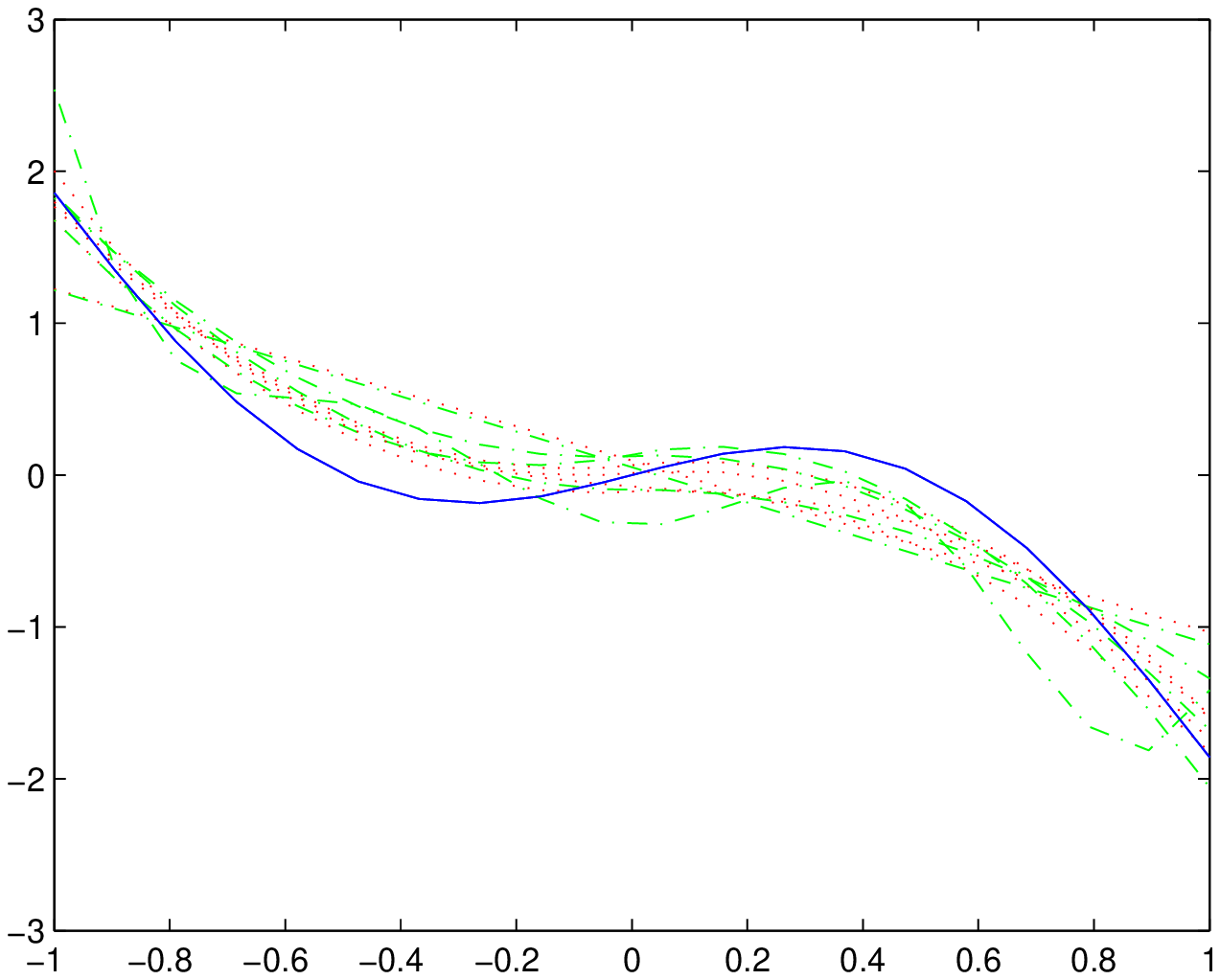}\\

\par\end{centering}

\centering{}-- : true function$\quad$-.-: first estimator$\quad\ldots$:
truncated estimator\\
$n=10^{4}$ et $\Delta=10^{-1}$
\end{figure}

\begin{figure}[p]
\caption{Model 4: Lévy process }

\begin{centering}
\[
b(x)=-2x,\;\sigma(x)=\xi(x)=1\textrm{ jumps Lévy }
\]
\label{Flo:figure_Model 4}
\par\end{centering}

\begin{centering}
\includegraphics[scale=0.5]{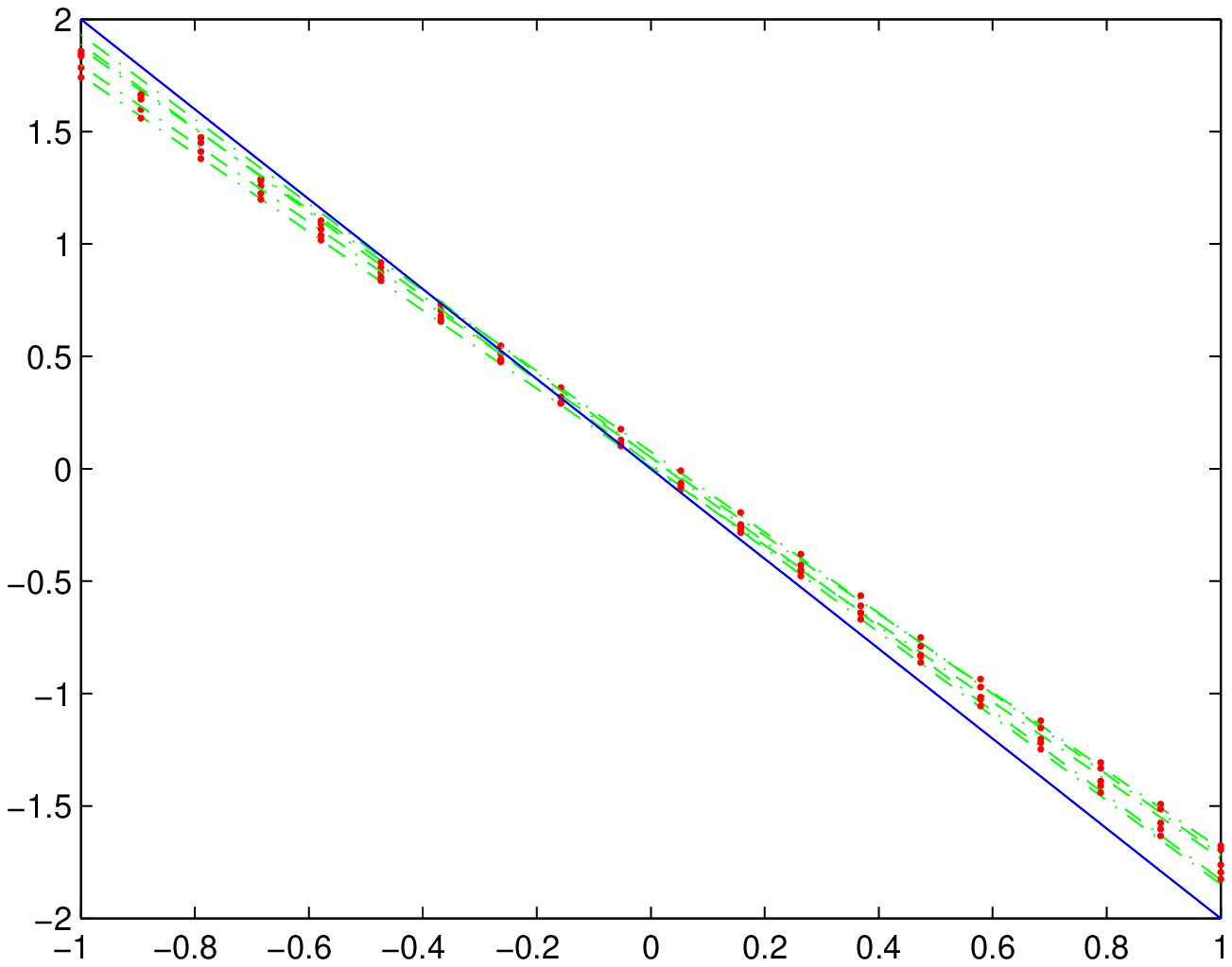}\\

\par\end{centering}

\centering{}-- : true function$\quad$-.-: first estimator$\quad\ldots$:
truncated estimator\\
$n=10^{4}$ et $\Delta=10^{-1}$
\end{figure}

In Figures \ref{Flo:figure_model 1}-\ref{Flo:figure_Model 4}, we
simulate 5 times the process $\left(X_{0},\ldots,X_{n\Delta}\right)$
for $\Delta=10^{-1}$ and $n=10^{4}$ and draw the obtained estimators.
The two adaptive estimators are nearly superposed, moreover, they
are close to the true function. 

In Tables \ref{Flo:table_Model1}-\ref{Flo:table_model4}, for each
value of $\left(n,\Delta\right)$, we simulate 50 trajectories of
$\left(X_{0},X_{\Delta},\ldots,X_{n\Delta}\right)$. For each path,
we construct the two adaptive estimators $\hat{b}_{\hat{m},\hat{r}}$
and $\tilde{b}_{\tilde{m},\tilde{r}}$ and we compute the empirical
errors: 
\[
err_{1}=\left\Vert \hat{b}_{\hat{m},\hat{r}}-b_{A}\right\Vert _{n}^{2}
\quad\textrm{and}\quad 
err_{2}=\left\Vert \tilde{b}_{\tilde{m},\tilde{r}}-b_{A}\right\Vert _{n}^{2}.
\]
In order to check that our algorithm is adaptive, we also compute
the minimal errors 
\[
emin_{1}=\min_{m,r}\left\Vert \hat{b}_{m,r}-b_{A}\right\Vert _{n}^{2}
\quad\textrm{and}\quad 
emin_{2}=\min_{m,r}\left\Vert \tilde{b}_{m,r}-b_{A}\right\Vert _{n}^{2}
\]
and the oracles $oracle_{i}=err_{i}/emin_{i}$. We give the means
$\hat{m}_{a}$, $\hat{r}_{a}$, $\tilde{m}_{a}$ and $\tilde{r}_{a}$
of the selected values $\hat{m}$, $\hat{r}$, $\tilde{m}$ and $\tilde{r}$.
The value $risk_{i}$ is the mean of $err_{i}$ over the 50 simulations
and $or_{i}$ is the mean of $oracle_{i}$. 
The computation time for one adaptive estimator varies from 0.1 second ($\Delta=10^{-1}$, $n=10^{3}$) to 30 seconds ($\Delta=10^{-1}$, $n=10^{4}$). 
The empirical risk is decreasing when the product $n\Delta$ is increasing,
which is coherent with the theoretical model. For Model 1, the two
estimators are equivalent. When the tails of $\nu$ become larger
(Models 2 and 3), the truncated estimator is better. The improvement
is also more significant when the discretization distance is smaller.
As on the first three models, the processes $L_{t}$ are compound Poisson
processes, these results were expected. The truncated estimator seems
also more robust: we do not observe aberrant values (like for the first
estimator in Table \ref{Flo:table_model2}). Those aberrant values may
be due to the fact that $b$ is not Lipschitz and then $b(X_{k\Delta})$
may be quite large, and to the non-exact simulation by an Euler scheme. 
For Model 4, the results are slightly better for the first estimator when $\Delta=0.1$, 
which is due to the fact that the remainder term is greater for the truncated estimator. 
When $\Delta=10^{-2}$, the risk of the truncated estimator is lower than for the first estimator.

\section{Proofs \label{sec:Proofs}}

Let us introduce the filtration 
\[
\rond{F}_{t}=\sigma\left(\eta,\left(W_{s}\right)_{0\leq s\leq t},\left(L_{s}\right)_{0\leq s\leq t}\right).
\]
The following result is very useful. It comes from \citet{dellacheriemeyer}
(Theorem 92 Chapter VII) and \citet{applebaum}, Theorem 4.4.23 p265
(Kunita's first inequality). 

\begin{result}[Burkholder-Davis-Gundy inequality]\label{res_burkholder}

We have that, for any $p\geq2$, 

\[
\E\left[\left.\sup_{s\in[t,t+h]}\left|\int_{t}^{s}\sigma(X_{u})dW_{u}\right|^{p}\right|\rond{F}_{t}\right]\leq C_{p}\left(\E\left[\left.\left|\int_{t}^{t+h}\sigma^{2}(X_{u})du\right|^{p/2}\right|\rond{F}_{t}\right]\right)\]
and, if $\int_{\mathbb{R}}\left|z\right|^{p}\nu(dz)<\infty$, as $\int_{\mathbb{R}}z^{2}\nu(dz)=1$:
\begin{eqnarray*}
\E\left[\left.\sup_{s\in[t,t+h]}\left|\int_{t}^{s}\xi(X_{u^{-}})dL_{u}\right|^{p}\right|\rond{F}_{t}\right] 
& \leq &
 C_{p}\E\left[\left.\left(\int_{t}^{t+h}\xi^{2}(X_{u})du\right)^{p/2}\right|\rond{F}_{\text{t}}\right]\\
 & + & 
C_{p}\E\left[\left.\left(\int_{t}^{t+h}\left|\xi(X_{u})\right|^{p}du\right)\right|\rond{F}_{t}\right]\int_{\mathbb{R}}\left|z\right|^{p}\nu(dz).
\end{eqnarray*}

\end{result}

% \subsection{Proof of Proposition \ref{prop_gloter}}
% 
% By Result \ref{res_burkholder}, there exists a constant $c_{p}$
% such that:
% \begin{eqnarray*}
% \E\left[\left.\sup_{s\in[t,t+h]}\left(X_{s}-X_{t}\right)^{2p}\right|\rond{F}_{t}\right] 
% & \leq & 
% c(p)\left(\left[\left.\left(\int_{t}^{t+h}\left|b(X_{s})\right|ds\right)^{2p}\right|\rond{F}_{t}\right]\right)\\
%  & + & 
% c(p)\E\left[\left.\left(\int_{t}^{t+h}\sigma^{2}(X_{s})ds\right)^{p}\right|\rond{F}_{t}\right]\\
%  & + &
%  c(p)\left(\E\left[\left.\left(\int_{t}^{t+h}\xi^{2}\left(X_{s}\right)ds\right)^{p}+\int_{t}^{t+h}\xi^{2p}(X_{s})ds\right|\rond{F}_{t}\right]\right).
% \end{eqnarray*}
% 
% Then, as $\xi$ and $\sigma$ are bounded and $b$ Lipschitz (and
% thus sub-linear), there exists a constant $C_{b}$ such that:
% \[
% \E\left[\left.\sup_{s\in[t,t+h]}\left(X_{s}-X_{t}\right)^{2p}\right|\rond{F}_{t}\right]
% \leq 
% c(p)\left(\sigma_{0}^{2p}h^{p}+\xi_{0}^{2p}(h+h^{p})\right)+c(p)h^{2p-1}C_{b}\int_{t}^{t+h}\E\left[\left.X_{s}^{2p}\right|\rond{F}_{t}\right]ds.
% \]
% As $(X_{t})$ is stationary, we obtain the expected result. 

\subsection{Proof of Theorem \ref{thme_risque_m_fixe} \label{sub:Proof-of-Theorem_1}}

By \eqref{eq:def_gamma} and \eqref{eq:def_risque}, we get: 
\begin{eqnarray*}
\gamma_{n}(t)=\frac{1}{n}\sum_{k=1}^{n}\left(Y_{k\Delta}-t(X_{k\Delta})\right)^{2} 
& = & 
\frac{1}{n}\sum_{k=1}^{n}\left(Y_{k\Delta}-b(X_{k\Delta})\right)^{2}+\left\Vert b-t\right\Vert _{n}^{2}\\
 & + & 
\frac{2}{n}\sum_{k=1}^{n}\left(Y_{k\Delta}-b(X_{k\Delta})\right)\left(b(X_{k\Delta})-t(X_{k\Delta})\right).
\end{eqnarray*}
 As, by definition, $\gamma_{n}(\hat{b}_{m})\leq\gamma_{n}(b_{m})$,
we obtain:
\[
\left\Vert \hat{b}_{m}-b\right\Vert _{n}^{2}
\leq
\left\Vert b_{m}-b\right\Vert _{n}^{2}+\frac{2}{n}\sum_{k=1}^{n}\left(Y_{k\Delta}-b(X_{k\Delta})\right)
\left(\hat{b}_{m}(X_{k\Delta})-b_{m}(X_{k\Delta})\right).
\]
By \eqref{eq:def_Y}, and as $\hat{b}_{m}$ and $b_{m}$ are supported
by $A$, 
\[
\left\Vert \hat{b}_{m}-b_{A}\right\Vert _{n}^{2}
\leq
\left\Vert b_{m}-b_{A}\right\Vert _{n}^{2}+\frac{2}{n}\sum_{k=1}^{n}\left(I_{k\Delta}+Z_{k\Delta}+T_{k\Delta}\right)
\left(\hat{b}_{m}(X_{k\Delta})-b_{m}(X_{k\Delta})\right).
\]
Let us introduce the unit ball 
\[
\rond{B}_{m}=\left\{ t\in S_{m},\;\left\Vert t\right\Vert _{\varpi}\leq1\right\} 
\quad\textrm{where}\quad
\left\Vert t\right\Vert _{\varpi}^{2}=\int_{A}t^{2}(x)\varpi(dx)
\]
and the englobing space 
$\rond{S}_{n}=\bigcup_{m\in\rond{M}_{n}}S_{m}$.
Let us consider the set 
\[
\Omega_{n}=\left\{ \omega,\;\forall t\in\rond{S}_{n}\;,\;
\left|\frac{\left\Vert t\right\Vert _{n}^{2}}{\left\Vert t\right\Vert _{\varpi}^{2}}-1\right|\leq\frac{1}{2}\right\} 
\]
where the norms $\left\Vert .\right\Vert _{\varpi}$ and $\left\Vert .\right\Vert _{n}$
are equivalent.

\paragraph{Step 1: bound of the risk on $\Omega_{n}$}

Thanks to the Cauchy-Schwartz inequality, we obtain that, on $\Omega_{n}$: 

\[
\left\Vert \hat{b}_{m}-b_{A}\right\Vert _{n}^{2}
\leq
\left\Vert b_{m}-b_{A}\right\Vert _{n}^{2}+\frac{1}{12}\left\Vert \hat{b}_{m}-b_{m}\right\Vert _{n}^{2}
+12\sum_{k=1}^{n}I_{k\Delta}^{2}+\frac{1}{12}\left\Vert \hat{b}_{m}-b_{m}\right\Vert _{\varpi}^{2}+12\sup_{t\in\rond{B}_{m}}\nu_{n}^{2}(t)
\]
 where 
\begin{equation}
\nu_{n}(t)=\frac{1}{n}\sum_{k=1}^{n}(Z_{k\Delta}+T_{k\Delta})t(X_{k\Delta}).
\label{eq:def_nun}
\end{equation}
On $\Omega_{n}$, by definition, we have: 
\[
\left\Vert \hat{b}_{m}-b_{m}\right\Vert _{n}^{2}
\leq2\left\Vert \hat{b}_{m}-b_{A}\right\Vert _{n}^{2}+2\left\Vert b_{m}-b_{A}\right\Vert _{n}^{2}
\quad\textrm{and}\quad
\left\Vert \hat{b}_{m}-b_{m}\right\Vert _{\varpi}^{2}\leq2\left\Vert \hat{b}_{m}-b_{m}\right\Vert _{n}^{2}.
\]
Thus we obtain: 
\[
\left\Vert \hat{b}_{m}-b_{A}\right\Vert _{n}^{2}
\leq3\left\Vert b_{m}-b_{A}\right\Vert _{n}^{2}+24\sum_{k=1}^{n}I_{k\Delta}^{2}+24\sup_{t\in\rond{B}_{m}}\nu_{n}^{2}(t).
\]
The following lemma is very useful. It is derived from Proposition \ref{prop_gloter} and Result \ref{res_burkholder}. 

\begin{lem}\label{lem_majoration_puissances}
\begin{enumerate}
\item $\E\left(I_{k\Delta}^{2}\right)\leq c\Delta$ and $\E\left(I_{k\Delta}^{4}\right)\leq c\Delta$.
\item $\E\left(\left.Z_{k\Delta}\right|\rond{F}_{k\Delta}\right)=0$, $\E\left(\left.Z_{k\Delta}^{2}\right|\rond{F}_{k\Delta}\right)\leq\sigma_{0}^{2}/\Delta$
and $\E\left(\left.Z_{k\Delta}^{4}\right|\rond{F}_{k\Delta}\right)\leq c/\Delta^{2}$.
\item $\E\left(\left.T_{k\Delta}\right|\rond{F}_{k\Delta}\right)=0$, $\E\left(\left.T_{k\Delta}^{2}\right|\rond{F}_{k\Delta}\right)\leq\xi_{0}^{2}/\Delta$
and $\E\left(\left.T_{k\Delta}^{4}\right|\rond{F}_{k\Delta}\right)\leq c/\Delta^{3}$.
\end{enumerate}
\end{lem}

By Lemma \ref{lem_majoration_puissances}, $\E\left[I_{k\Delta}^{2}\right]\leq\Delta$.
It remains to bound
 $\E\left[\sup_{t\in\rond{B}_{m}}\nu_{n}^{2}(t)\right].$
We consider an orthonormal basis 
$\left(\varphi_{\lambda}\right)_{\lambda\in\Lambda_{m}}$
of $S_{m}$ for the $L_{\varpi}^{2}$-norm with $\vert\Lambda_{m}\vert=D_{m}$.
Any function $t\in S_{m}$ can be written 
$t=\sum_{\lambda\in\Lambda_{m}}a_{\lambda}\varphi_{\lambda}$
and 
$\left\Vert t\right\Vert _{\varpi}^{\text{2}}=\sum_{\lambda\in\Lambda_{m}}a_{\lambda}^{2}.$
Then: 
\begin{eqnarray*}
\sup_{t\in\rond{B}_{m}}\nu_{n}^{2}(t) 
& = & 
\sup_{\sum_{\lambda}a_{\lambda}^{2}\leq1}\left(\sum_{\lambda\in\Lambda_{m}}a_{\lambda}\nu_{n}\left(\varphi_{\lambda}\right)\right)^{2}\\
 & \leq & 
\sup_{\sum_{\lambda}a_{\lambda}^{2}\leq1}\left(\sum_{\lambda\in\Lambda_{m}}a_{\lambda}^{2}\right)
\left(\sum_{\lambda\in\Lambda_{m}}\nu_{n}^{2}\left(\varphi_{\lambda}\right)\right)\\
 & = & 
\sum_{\lambda\in\Lambda_{m}}\nu_{n}^{2}\left(\varphi_{\lambda}\right).
\end{eqnarray*}
It remains to bound 
$\E\left(\nu_{n}^{2}\left(\varphi_{\lambda}\right)\right)$.
By \eqref{eq:def_nun}, 
\begin{eqnarray*}
\E\left[\nu_{n}^{2}(\varphi_{\lambda})\right] 
& = & 
\frac{1}{n^{2}}\sum_{k=1}^{n}\E\left[\varphi_{\lambda}^{2}(X_{k\Delta})\E\left[\left.(Z_{k\Delta}+T_{k\Delta})^{2}\right|\rond{F}_{k\Delta}\right]\right]\\
 & + &
 \frac{2}{n^{2}}\sum_{k<l}^{n}\E\left[(Z_{k\Delta}+T_{k\Delta})\varphi_{\lambda}(X_{k\Delta})\varphi_{\lambda}(X_{l\Delta})
\E\left[\left.Z_{l\Delta}+T_{l\Delta}\right|\rond{F}_{l\Delta}\right]\right]
\end{eqnarray*}
Thanks to Lemma \ref{lem_majoration_puissances}, the second term
of this inequality is null and we obtain, as $\int_{\mathbb{R}}\varphi_{\lambda}^{2}(x)\varpi(dx)=1$:
\[
\E\left[\nu_{n}^{2}(\varphi_{\lambda})\right]\leq
\frac{2(\sigma_{0}^{2}+\xi_{0}^{2})}{n^{2}\Delta}\sum_{k=1}^{n}\E\left[\varphi_{\lambda}^{2}(X_{k\Delta})\right]
=\frac{2(\sigma_{0}^{2}+\xi_{0}^{2})}{n\Delta}.
\]
Therefore: 
 \[
\E\left[\left\Vert \hat{b}_{m}-b_{A}\right\Vert _{n}^{2}\units{\Omega_{n}}\right]
\leq3\left\Vert b_{m}-b_{A}\right\Vert _{n}^{2}+48(\sigma_{0}^{2}+\xi_{0}^{2})\frac{D_{m}}{n\Delta}+C\Delta.
\]

% \begin{req}
%  If $\sigma$ and $\xi$ are not bounded, our estimator converges with the same rate. Indeed, we can bound $\E\left(\sup_{t\in \rond{B}_m}\nu_n^2(t)\right)$ 
% in a slightly different way. Let us consider an orthornormal basis $\left\{\varphi_{\lambda}\right\}_{\lambda\in\Lambda_{m}}$ of $S_m$. Any $t$ beloging to 
% $\rond{B}_m$ can be written
% \[t=\sum_{\lambda\in\Lambda_m}a_\lambda \varphi_{\lambda}\quad\textrm{ with}\quad  
%  \sum_{\lambda\in\Lambda_m}a_{\lambda}^2=\norm{t}_{L^2}^2\leq \pi_1^{-1}\norm{t}_{  \varpi}^2\leq \pi_1^{-1}
% \]
% Then, by Assumption A\ref{hypo_espaces} \it{3}, 
% \begin{eqnarray*}
%  \E\left(\sup_{t\in\rond{B}_m}\nu_n^2(t)\right)&\leq &\frac{\pi_1^{-1}}{n^2}\sum_{k=1}^{n}\E\left(\sum_{\lambda\in\Lambda_m}\varphi_{\lambda}^2(X_{k\Delta})
%  \E\left((Z_{k\Delta}+T_{k\Delta})^2\left|\rond{F}_{k\Delta}\right.\right)\right)\\
% &\leq& \frac{\phi_1 D_m}{n\Delta}\E\left(\sigma^2(X_{0})+\xi^2(X_{0})\right)
% \end{eqnarray*}
% 
% 
% 
% \end{req}

\paragraph{Step 2: bound of the risk on $\Omega_{n}^{c}.$ }

The process $\left(X_{t}\right)_{t\geq0}$ is exponentially $\beta$-mixing,
$\pi$ is bounded from below and above and $n\Delta\rightarrow\infty$.
The following result is proved for $\xi=0$ for instance in \citet{comtegenon2007} for diffusion processes,
but as it relies only on the $\beta$-mixing property, we can apply
it. 

\begin{result}\label{res_borne_omega}

\[
\mathbb{P}\left[\Omega_{n}^{c}\right]\leq\frac{1}{n^{3}}.
\]

\end{result}

Let us set $e=\left(e_{\Delta},\ldots,e_{n\Delta}\right)^{*}$ where
$e_{k\Delta}:=Y_{k\Delta}-b(X_{k\Delta})=I_{k\Delta}+Z_{k\Delta}+T_{k\Delta}$
and 
$\Pi_{m}Y=\Pi_{m}\left(Y_{\Delta},\ldots,Y_{n\Delta}\right)^{*}=\left(\hat{b}_{m}(X_{0}),\ldots,\hat{b}_{m}(X_{n\Delta})\right)^{*}$
where $\Pi_{m}$ is the Euclidean orthogonal projection over $S_{m}$.
Then 
\begin{eqnarray*}
\left\Vert \hat{b}_{m}-b_{A}\right\Vert _{n}^{2} 
& = & 
\left\Vert \Pi_{m}Y-b_{A}\right\Vert _{n}^{2}=\left\Vert \Pi_{m}b_{A}-b_{A}\right\Vert _{n}^{2}+\left\Vert \Pi_{m}Y-\Pi_{m}b_{A}\right\Vert _{n}^{2}\\
 & \leq & 
\left\Vert b_{A}\right\Vert _{n}^{2}+\left\Vert e\right\Vert _{n}^{2}.
\end{eqnarray*}
According to Lemma \ref{lem_majoration_puissances}, Result \ref{res_borne_omega}
and the Cauchy-Schwarz inequality,
 \[
\E\left[\left\Vert e\right\Vert _{n}^{2}\units{\Omega_{n}^{c}}\right]
\leq
\left(\E\left[\left\Vert e\right\Vert _{n}^{4}\right]\right)^{1/2}
\left(\mathbb{P}\left(\Omega_{n}^{c}\right)\right)^{1/2}\leq\frac{C}{\left(\Delta^{3}n^{3}\right)^{1/2}}\leq\frac{C}{n\Delta}
\]
and, as $b$ is bounded on the compact set $A$, 
\[
\E\left[\left\Vert b_{A}\right\Vert _{n}^{2}\units{\Omega_{n}^{c}}\right]
\leq
\left(\E\left[\left\Vert b_{A}\right\Vert _{n}^{4}\right]\mathbb{P}\left(\Omega_{n}^{c}\right)\right)^{1/2}\lesssim\frac{1}{n^{3/2}}.
\]
Collecting the results, we get: 
\[
\E\left[\left\Vert \hat{b}_{m}-b_{A}\right\Vert _{n}^{2}\units{\Omega_{n}^{c}}\right]\lesssim\frac{1}{n\Delta}
\]
which ends the proof of Theorem \ref{thme_risque_m_fixe}.

\subsection{Proof of Theorem \ref{thme_adaptatif} \label{sub:Proof-of-Theorem_adaptatif}}

The bound of the risk on $\Omega_{n}^{c}$ is done exactly in the
same way as for the non adaptive estimator. It remains thus to bound
the risk on $\Omega_{n}$. As in the previous proof, we get: 
\begin{eqnarray*}
\left\Vert \hat{b}_{\hat{m}}-b_{A}\right\Vert _{n}^{2}\units{\Omega_{n}} 
& \leq & 
3\left\Vert b_{m}-b_{A}\right\Vert _{n}^{2}+\frac{24}{n}\sum_{k=1}^{n}I_{k\Delta}^{2}+2pen(m)-2pen(\hat{m})\\
 & + & 24\sup_{t\in\rond{B}_{m,\hat{m}}}\nu_{n}^{2}(t)
\end{eqnarray*}

where $\rond{B}_{m,m'}$ is the unit ball (for the $L_{\varpi}^{2}$-norm)
of the subspace $S_{m}+S_{m'}$: $\rond{B}_{m,m'}=\left\{ t\in S_{m}+S_{m'},\:\left\Vert t\right\Vert _{\varpi}\leq1\right\} $.
Let us introduce a function $p(m,m')$ such that $12p(m,m')=pen(m)+pen(m')$.
We obtain that, on $\Omega_{n}$, for any $m\in\rond{M}_{n}$: 
\begin{eqnarray*}
\left\Vert \hat{b}_{\hat{m}}-b_{A}\right\Vert _{n}^{2} 
& \leq &
 3\left\Vert b_{m}-b_{A}\right\Vert _{n}^{2}+\frac{24}{n}\sum_{k=1}^{n}I_{k\Delta}^{2}+4pen(m)\\
 & + & 24\sup_{t\in\rond{B}_{m,\hat{m}}}\left(\nu_{n}^{2}(t)-p(m,\hat{m})\right).
\end{eqnarray*}
It remains to bound 
\[
\E\left[\sup_{t\in\rond{B}_{m,\hat{m}}}\nu_{n}^{2}(t)-p(m,\hat{m})\right]\leq\sum_{m'}\E\left[\sup_{t\in\rond{B}_{m,m'}}\nu_{n}^{2}(t)-p(m,m')\right]_{+}.
\]
For this purpose, we use the following proposition proved in \citet{applebaum}
(Corollary 5.2.2 ). 

\begin{prop}[exponential martingale]\label{prop_martingale_exponentielle}

Let $(Y_{t})_{t\geq0}$ satisfy: 
\[
Y_{t}=\int_{0}^{t}F_{s}dW_{s}+\int_{0}^{t}K_{s}dL_{s}-\int_{0}^{t}\left[\frac{F_{s}^{2}}{2}+\int_{\mathbb{R}}\left(e^{K_{s}z}-1-K_{s}z\right)\nu(dz)\right]ds
\]
where $F_{s}$ and $K_{s}$ are locally integrable and predictable
processes. If for any $t>0$, 
\[
\E\left[\int_{0}^{t}\int_{\vert z\vert>1}\left|e^{K_{s}z}-1\right|\nu(dz)ds\right]<\infty,
\]
 then $e^{Y_{t}}$ is a $\rond{G}_{t}$-local martingale where 
$\rond{G}_{t}=\sigma(W_{s},L_{s},0\leq s\leq t)$.

\end{prop}

For any 
$\ep\leq\ep_{1}:=(\lambda\wedge1)/(2\left\Vert t\right\Vert _{\infty}\xi_{0})$
where $\lambda$ is defined in Assumption A\ref{hypo_loi_xi}, for
any $t\geq0$ 
\[
\int_{0}^{t}\int_{\left|z\right|\geq1}\left(\exp(\ep t(X_{k\Delta})\xi(X_{s})z)-1\right)\nu(dz)\units{s\in]k\Delta,(k+1)\Delta]}ds<\infty.
\]
 Let us introduce the two Markov processes 
\[
A_{\ep,t}:=\ep^{2}\sum_{k=0}^{n}t^{2}(X_{k\Delta})\int_{0}^{t}\sigma^{2}(X_{s})\units{s\in]k\Delta,(k+1)\Delta]}ds
\]
and 

\[
B_{\ep,t}:=\sum_{k=0}^{n}\int_{0}^{t}\int_{\mathbb{R}}\left(\exp\left(\ep t(X_{k\Delta})\xi(X_{s})z\right)
-\ep t(X_{k\Delta})\xi(X_{s})z-1\right)\units{s\in]k\Delta,(k+1)\Delta]}\nu(dz)ds
\]
and the following martingale: 

\[
M_{t}=\int_{0}^{t}\sum_{k=0}^{n}\units{s\in]k\Delta,(k+1)\Delta]}t(X_{k\Delta^{-}})\left(\sigma(X_{s})dW_{s}+\xi(X_{s^{-}})dL_{s}\right).\]
By Proposition \ref{prop_martingale_exponentielle}, 
\[
Y_{\ep,s}:=\ep M_{s}-A_{\ep,s}-B_{\ep,s}
\]
 is such that $e^{Y_{\ep,s}}$ is a local martingale.

\paragraph{Bound of $A_{\ep,s}$ and $B_{\ep,s}$.}

We obtain easily that $A_{\ep,s}\leq A_{\ep,(n+1)\Delta}\leq\ep^{2}n\Delta\left\Vert t\right\Vert _{n}^{2}\sigma_{0}^{2}$.
Under Assumption A\ref{hypo_loi_xi}, $\xi$ is constant or $\nu$
is symmetric, and therefore 
\[
B_{\ep,s}\leq B_{\ep,(n+1)\Delta}
\leq\Delta\sum_{k=0}^{n}\int_{\mathbb{R}}\left(\exp\left(\ep t(X_{k\Delta})\xi_{0}z\right)-\ep t(X_{k\Delta^{-}})\xi_{0}z-1\right)\nu(dz).
\]
As $\int_{\mathbb{R}}z^{2}\nu(dz)=1$, for any $\alpha\leq1$, 
\[
\int_{-1}^{1}\left(\exp\left(\alpha z\right)-\alpha z-1\right)\nu(dz)\leq\alpha^{2}\int_{-1}^{1}z^{2}\nu(dz)\leq\alpha^{2}.
\]
Moreover, by integration by parts, for any $\alpha\leq(1\wedge\lambda)/2$,
\begin{eqnarray*}
\int_{[-1,1]^{c}}\left(\exp\left(\alpha z\right)-\alpha z-1\right)\nu(dz) 
& \leq & 
\left(e^{\alpha}-\alpha-1\right)\nu([1,+\infty[)+\left(e^{-\alpha}+\alpha-1\right)\nu(]-\infty,-1])\\
& + & 
\int_{1}^{+\infty}\alpha\left(e^{\alpha z}-1\right)\nu([-z,z]^{c})dz\\
\end{eqnarray*}
By assumption A\ref{hypo_loi_xi}, $\nu([-z,z]^c)\leq Ce^{-\lambda z}$ and then 
\[
 \int_{[-1,1]^{c}}\left(\exp\left(\alpha z\right)-\alpha z-1\right)\nu(dz) 
  \leq  2\alpha^{2}\nu\left([-1,1]^{c}\right)+Ce^{-\lambda}\frac{\alpha}{\lambda}\left(\frac{e^{\alpha}}{1-\alpha/\lambda}-1\right)\leq C'\alpha^{2}.
\]

Then 
$B_{\ep,s}\lesssim n\Delta\varepsilon^{2}\xi_{0}^{2}\left\Vert t\right\Vert _{n}^{2}$.
There exists a constant $c$ such that, for any $\ep<\ep_{1}$, 
\[
A_{\ep,s}+B_{\ep,s}\leq c\frac{n\Delta\ep^{2}\left(\sigma_{0}^{2}+\xi_{0}^{2}\right)\left\Vert t\right\Vert _{n}^{2}}{\left(1-\ep/\ep_{1}\right)}.
\]

\paragraph{Bound of $\mathbb{P}\left(\nu_{n}(t)\geq\eta,\:\left\Vert t\right\Vert _{n}^{2}\leq\zeta^{2}\right)$.}

The process $\exp(Y_{\ep,t})$ is a local martingale, then there exists
an increasing sequence $(\tau_{N})$ of stopping times such that 
$\lim_{N\rightarrow\infty}\tau_{N}=\infty$
and 
$\exp(Y_{\ep,t\wedge\tau_{N}})$
 is a $\rond{F}_{t}$-martingale.
For any $\ep<\ep_{1}$, and all $N$, 
\begin{eqnarray*}
E & := & 
\mathbb{P}\left(M_{(n+1)\Delta\wedge\tau_{N}}\geq n\Delta\eta,\:\left\Vert t\right\Vert _{n}^{2}\leq\zeta^{2}\right)\\
 & \leq & 
\mathbb{P}\left(M_{(n+1)\Delta\wedge\tau_{N}}\geq n\Delta\eta,\:
 A_{(n+1)\Delta\wedge\tau_{N}}+B_{(n+1)\Delta\wedge\tau_{N}}
\leq\frac{cn\Delta\ep^{2}\left(\sigma_{0}^{2}+\xi_{0}^{2}\right)\zeta^{2}}{\left(1-\ep/\ep_{1}\right)}\right)\\
 & \leq & 
\E\left(\exp(Y_{\ep,(n+1)\Delta\wedge\tau_{N}})\right)\exp\left(-n\Delta\eta\ep
+\frac{cn\Delta\ep^{2}\left(\xi_{0}^{2}+\sigma_{0}^{2}\right)\zeta^{2}}{\left(1-\ep/\ep_{1}\right)}\right).
\end{eqnarray*}
As $\exp(Y_{\ep,t\wedge\tau_{N}})$ is a martingale, $\E\left(\exp(Y_{\ep,t\wedge\tau_{N}})\right)=1$
and 
\[
E\leq\exp\left(-n\Delta\eta\ep+\frac{cn\Delta\ep^{2}\left(\xi_{0}^{2}+\sigma_{0}^{2}\right)\zeta^{2}}{\left(1-\ep/\ep_{1}\right)}\right).
\]
Letting $N$ tend to infinity, by dominated convergence, and as $\nu_{n}(t)=n\Delta M_{(n+1)\Delta}$,
we obtain that 
\[
\mathbb{P}\left(\nu_{n}(t)\geq\eta,\:
\left\Vert t\right\Vert _{n}^{2}\leq\zeta^{2}\right)
\leq\exp\left(-n\Delta\eta\ep+\frac{cn\Delta\ep^{2}\left(\xi_{0}^{2}+\sigma_{0}^{2}\right)\zeta^{2}}{\left(1-\ep/\ep_{1}\right)}\right).
\]
It remains to minimise this inequality in $\varepsilon$. Let us set
\[
\text{\ensuremath{\varepsilon}}=\frac{\eta}{2c\left(\sigma_{0}^{2}+\xi_{0}^{2}\right)\zeta^{2}/\Delta+\eta/\varepsilon_{1}}<\varepsilon_{1}.
\]
We get: 
\[
\mathbb{P}\left(\nu_{n}(t)\geq\eta,\;
\left\Vert t\right\Vert _{n}^{2}\leq\zeta^{2}\right)
\leq\exp\left(-\frac{\eta^{2}n\Delta}{4c\left(\left(\sigma_{0}^{2}+\xi_{0}^{2}\right)\zeta^{2}+c'\eta\xi_{0}\left\Vert t\right\Vert _{\infty}\right)}\right).
\]

The following lemma concludes the proof. It is proved thanks to a
$L_{\varpi}^{2}-L^{\infty}$ chaining technique. See \citet{comte2001},
proof of Proposition 4, and \citet{schmisser_these}, Appendix D.3. 

\begin{lem}

There exists a constant $\kappa$ such that: 

\[
\E\left[\sup_{t\in\rond{B}_{m,m'}}\nu_{n}^{2}(t)-p(m,m')\right]\lesssim\kappa(\xi_{0}^{2}+\sigma_{0}^{2})\frac{D^{3/2}}{n\Delta}e^{-\sqrt{D}}
\]
where $D=\dim(S_{m}+S_{m'})$. 

\end{lem}

As $\sum_{D}D^{3/2}e^{-\sqrt{D}}\leq\sum_{k=0}^{+\infty}k^{3/2}e^{-\sqrt{k}}<\infty$,
we obtain that 
\[
\E\left[\sup_{t\in\rond{B}_{m,\hat{m}}}\nu_{n}^{2}(t)-p(m,\hat{m})\right]
\leq
\sum_{m'\in\rond{M}_{n}}\E\left[\sup_{t\in\rond{B}_{m,m'}}\nu_{n}^{2}(t)-p(m,m')\right]\lesssim\kappa\frac{\xi_{0}^{2}+\sigma_{0}^{2}}{n\Delta}.
\]

\subsection{Proof of Theorem \ref{thme_sauts_coupes_m_fixe } }

We recall that 
\[
\Omega_{X,k}=\left\{ \omega,\;\left|X_{(k+1)\Delta}-X_{k\Delta}\right|\leq C_{\Delta}
=\left(b_{max}+3\right)\Delta+\left(\sigma_{0}+4\xi_{0}\right)\sqrt{\Delta}\ln(n)\right\} .
\]
Let us introduce the set 

\[
\Omega_{N,k}=\left\{ \omega,\: N_{k\Delta}^{'}=0\right\}
 \]
where $N_{k\Delta}^{'}$ is the number of jumps of size larger than
$\Delta^{1/4}$ occurring in the time interval $]k\Delta,(k+1)\Delta]$:
\[
N'_{k\Delta}=\mu\left(]k\Delta,(k+1)\Delta]\;,\:\left[-\Delta^{1/4},\Delta^{1/4}\right]^{c}\right).
\]
We have that 
\begin{eqnarray*}
\tilde{Y}_{k\Delta} & = & Y_{k\Delta}\units{\Omega_{X,k}}\units{X_{k\Delta}\in A}\\
 & = & b_{A}(X_{k\Delta})-b_{A}(X_{k\Delta})\units{\Omega_{X,k}^{c}\cap(X_{k\Delta}\in A)}
+I_{k\Delta}\units{\Omega_{X,k}\cap(X_{k\Delta}\in A)}+\tilde{Z}_{k\Delta}+\tilde{T}_{k\Delta}\\
 & + & \left(Z_{k\Delta}+T_{k\Delta}\right)\units{\Omega_{X,k}\cap\Omega_{N,k}^{c}\cap(X_{k\Delta}\in A)}
+\E\left(\left.\left(Z_{k\Delta}+T_{k\Delta}\right)\units{\Omega_{X,k}\cap\Omega_{N,k}\cap(X_{k\Delta}\in A)}\right|\rond{F}_{k\Delta}\right).
\end{eqnarray*}
where 
\[
\tilde{Z}_{k\Delta}=Z_{k\Delta}\units{\Omega_{X,k}\cap\Omega_{N,k}\cap(X_{k\Delta}\in A)}
-\E\left(\left.Z_{k\Delta}\units{\Omega_{X,k}\cap\Omega_{N,k}\cap(X_{k\Delta}\in A)}\right|\rond{F}_{k\Delta}\right)
\]
 and 
\[
\tilde{T}_{k\Delta}=T_{k\Delta}\units{\Omega_{X,k}\cap\Omega_{N,k}\cap(X_{k\Delta}\in A)}
-\E\left(\left.T_{k\Delta}\units{\Omega_{X,k}\cap\Omega_{N,k}\cap(X_{k\Delta}\in A)}\right|\rond{F}_{k\Delta}\right).
\]
As previously, we only bound the risk on $\Omega_{n}$. Let us set
\[
\tilde{\nu}_{n}(t):=\frac{1}{n}\sum_{k=1}^{n}t(X_{k\Delta})\left(\tilde{Z}_{k\Delta}+\tilde{T}_{k\Delta}\right).
\]
We have that 

\begin{eqnarray*}
\left\Vert \tilde{b}_{m}-b_{A}\right\Vert _{n}^{2}\units{\Omega_{n}}
 & \leq & 
3\left\Vert b_{m}-b_{A}\right\Vert _{n}^{2}+24\sup_{t\in\rond{B}_{m}}\tilde{\nu}_{n}^{2}(t)
+\frac{224}{n}\sum_{k=1}^{n}\left(I_{k\Delta}^{2}+b_{A}^{2}(X_{k\Delta})\units{\Omega_{X,k}^{c}}\right)\\
 & + & \frac{224}{n}\sum_{k=1}^{n}\left(Z_{k\Delta}^{2}+T_{k\Delta}^{2}\right)\units{\Omega_{X,k}\cap\Omega_{N,k}^{c}\cap(X_{k\Delta}\in A)}\\
 & + & \frac{224}{n}\sum_{k=1}^{n}\left(\E\left[\left.\left(Z_{k\Delta}+T_{k\Delta}\right)\units{\Omega_{X,k}\cap\Omega_{N,k}\cap(X_{k\Delta}\in A)}\right|
\rond{F}_{k\Delta}\right]\right)^{2}.
\end{eqnarray*}
The following lemma is proved later. 

\begin{lem}\label{lem_majoration_puissances_sauts_coupes}
\begin{enumerate}
\item $\mathbb{P}(\Omega_{X,k}^{c}\cap(X_{k\Delta}\in A))\lesssim\Delta^{1-\beta/2}$.
\item $\mathbb{P}(\Omega_{X,k}\cap\Omega_{N,k}^{c}\cap(X_{k\Delta}\in A))\lesssim\Delta^{2-\beta/2}$. 
\item $\left(\E\left[\left.\left(Z_{k\Delta}+T_{k\Delta}\right)\units{\Omega_{N,k}\cap\Omega_{X,k}\cap(X_{k\Delta}\in A)}\right|\rond{F}_{k\Delta}\right]\right)^{2}\lesssim\ln^{2}(n)\Delta^{1-\beta/2}$. 
\end{enumerate}
\end{lem}

According to Lemma \ref{lem_majoration_puissances}, $\E(I_{k\Delta}^{2})\leq\Delta$.
As $b$ is bounded on the compact set $A$, 
$\E\left[b_{A}^{2}(X_{k\Delta})\units{\Omega_{X,k}^{c}}\right]\lesssim\mathbb{P}(\Omega_{X,k}^{c})\lesssim\Delta^{1-\beta/2}.$
Moreover, on $\Omega_{X,k}$, 

\begin{eqnarray*}
\left(Z_{k\Delta}+T_{k\Delta}\right)^{2}\units{\Omega_{X,k}\cap(X_{k\Delta}\in A)}
 & = & 
\left(\frac{X_{(k+1)\Delta}-X_{k\Delta}}{\Delta}-b_{A}(X_{k\Delta})-I_{k\Delta}\right)^{2}\units{\Omega_{X,k}}\units{X_{k\Delta}\in A}\\
 & \lesssim & \frac{\ln^{2}(n)}{\Delta}+b_{A}^{2}(X_{k\Delta})+I_{k\Delta}^{2}
\end{eqnarray*}
and then
 \begin{eqnarray*}
E & := & \E\left[\left(Z_{k\Delta}+T_{k\Delta}\right)^{2}\units{\Omega_{X,k}\cap\Omega_{N,k}^{c}\cap(X_{k\Delta}\in A)}\right]\\
 & \lesssim & \left(\frac{\ln^{2}(n)}{\Delta}+b_{max}^{2}\right)\mathbb{P}\left(\Omega_{X,k}\cap\Omega_{N,k}^{c}\cap(X_{k\Delta}\in A)\right)
+\E\left(I_{k\Delta}^{2}\right)\\
 & \lesssim & \ln^{2}(n)\Delta^{1-\beta/2}.
\end{eqnarray*}
It remains to bound 
$\E\left(\sup_{t\in\rond{B}_{m}}\tilde{\nu}_{n}^{2}(t)\right)$.
In the same way as in Subsection \ref{sub:Proof-of-Theorem_1}, we
get: 
\begin{eqnarray*}
\E\left(\sup_{t\in\rond{B}_{m}}\tilde{\nu}_{n}^{2}(t)\right) 
& \leq & 
\sum_{\lambda\in\Lambda_{m}}\E\left(\tilde{\nu}_{n}^{2}(\varphi_{\lambda})\right)
\leq\frac{2D_{m}}{n}\E\left(\tilde{Z}_{\Delta}^{2}+\tilde{T}_{\Delta}^{2}\right)\\
 & \leq & \frac{2D_{m}}{n}\E\left(Z_{\Delta}^{2}+T_{\Delta}^{2}\right)\leq2\left(\sigma_{0}^{2}+\xi_{0}^{2}\right)\frac{D_{m}}{n\Delta}.
\end{eqnarray*}
We have that $\E\left(\tilde{Z}_{\Delta}^2\right)\leq \E\left(Z^2_{\Delta}\right)\leq \frac{\sigma_0^2}{\Delta}$. Moreover, 
\begin{eqnarray*}
 \E\left(\tilde{T}_{k\Delta}^2\right)&\lesssim& \E\left(T_{k\Delta}^2\units{\Omega_{X,k}\cap\Omega_{N,k}}\right)
-\left(\E\left(T_{k\Delta}\units{\Omega_{X,k}\cap\Omega_{N,k}}\right)\right)^2\\
&\lesssim& \E\left(T_{k\Delta}^2\units{\Omega_{N,k}}\right)+\ln^2(n)\Delta^{1-\beta/2}\\
&\lesssim& \Delta^{1/2-\beta/4}. 
\end{eqnarray*}
Then $\E\left(\sup_{t\in\rond{B}_{m}}\tilde{\nu}_{n}^{2}(t)\right) \leq (n\Delta)^{-1}D_m (\sigma_0^2+o(1))$.

\subsubsection{Proof of Lemma \ref{lem_majoration_puissances_sauts_coupes}}

\begin{result}\label{result_blumenthal}

Let $\beta$ be the Blumenthal-Getoor index of $L_{t}$. Then: 
\[
\nu([-z,z]^{c})\lesssim z^{-\beta}\quad,\quad
\int_{\left|x\right|\leq z\wedge a_{0}}x^{2}\nu(dx)\lesssim z^{2-\beta}
\quad\textrm{and}\quad
\int_{\left|x\right|\leq z\wedge a_{0}}x^{4}\nu(dx)\lesssim z^{4-\beta}.
\]
The constant $a_{0}$ is defined in A\ref{hypo_nu_reguliere}. 

\end{result}

\paragraph{Bound of $\mathbb{P}(\Omega_{X,k}^{c}\cap(X_{k\Delta}\in A))$.}

We have: 
\[
\mathbb{P}\left(\Omega_{X,k}^{c}\cap(X_{k\Delta}\in A)\right)
=\mathbb{P}\left(\left\{ \left|X_{(k+1)\Delta}-X_{k\Delta}\right|>C_{\Delta}\right\} \cap(X_{k\Delta}\in A)\right).
\]
We know that $X_{(k+1)\Delta}-X_{k\Delta}=b(X_{k\Delta})+I_{k\Delta}+Z_{k\Delta}+T_{k\Delta}$.
Then \begin{eqnarray*}
\mathbb{P}\left(\Omega_{X,k}^{c}\cap(X_{k\Delta}\in A)\right)&\leq&\mathbb{P}\left(\left|\Delta I_{k\Delta}\right|\geq\Delta\right)\\
&+&\mathbb{P}\left(\left|\Delta Z_{k\Delta}\right|\geq\sigma_{0}\sqrt{\Delta}\ln(n)\right)
+\mathbb{P}\left(\left|\Delta T_{k\Delta}\right|\geq\xi_{0}\sqrt{\Delta}\ln(n)\right).
\label{eq:decomposition_P_omega_X}
\end{eqnarray*}
By a Markov inequality and Lemma \ref{lem_majoration_puissances},
we obtain: 
\begin{equation}
\mathbb{P}\left(\left|\Delta I_{k\Delta}\right|\geq\Delta\right)\leq\frac{\E\left(\Delta^{2}I_{k\Delta}^{2}\right)}{\Delta^{2}}\lesssim\Delta.
\label{eq:majoration_I}
\end{equation}
By Proposition \ref{prop_martingale_exponentielle}, the process 
$\exp\left(c\int_{0}^{t}\sigma(X_{s^{-}})dW_{s}-c^{2}\int_{0}^{t}\sigma^{2}(X_{s})ds\right)$
is a local martingale (as $\sigma$ is bounded, it is in fact a martingale,
see \citet{liptser_shiryaev}, pp 229-232). Then, by a Markov inequality:
\begin{equation}
\mathbb{P}\left(\left|\Delta Z_{k\Delta}\right|\geq\sigma_{0}\sqrt{\Delta}\ln(n)\right)
\leq\frac{2}{n}\E\left[\exp\left(\frac{\sqrt{\Delta}Z_{k\Delta}}{\sigma_{0}}\right)\right]\lesssim\frac{1}{n}.
\label{eq:majoration_Z}
\end{equation}
To bound inequality \eqref{eq:decomposition_P_omega_X}, it remains
to bound 
$\mathbb{P}\left(\left|\Delta T_{k\Delta}\right|\geq\xi_{0}\sqrt{\Delta\ln(n)}\right)$.
Let us set 
\[
T_{k\Delta}=T_{k\Delta}^{(1)}+T_{k\Delta}^{(2)}+T_{k\Delta}^{(3)}
\quad\textrm{where}\quad
 T_{k\Delta}^{(i)}=\frac{1}{\Delta}\int_{k\Delta}^{(k+1)\Delta}\xi(X_{s^{-}})dL_{s}^{(i)}
\]
with \begin{gather*}
L_{t}^{(1)}=\int_{0}^{t}\int_{[-\sqrt{\Delta},\sqrt{\Delta}]}z\tilde{\mu}(ds,dz),
\quad
 L_{t}^{(2)}=\int_{0}^{t}\int_{[-\Delta^{1/4},-\sqrt{\Delta}]\cup[\sqrt{\Delta},\Delta^{1/4}]}z\tilde{\mu}(ds,dz)\\
L_{t}^{(3)}=\int_{0}^{t}\int_{[-\Delta^{1/4},\Delta^{1/4}]^{c}}z\tilde{\mu}(ds,dz).
\end{gather*}
Let us set 
$N_{k\Delta}^{''}=\mu\left(]k\Delta,(k+1)\Delta],\left[-\sqrt{\Delta},\sqrt{\Delta}\right]^{c}\right)$.
By Result \ref{result_blumenthal}, we have: 
\[
\mathbb{P}\left(\left|T_{k\Delta}^{(2)}+T_{k\Delta}^{(3)}\right|>0\right)
=\mathbb{P}\left(N_{k\Delta}^{''}\geq1\right)\lesssim\Delta\nu\left(\left[-\sqrt{\Delta},\sqrt{\Delta}\right]^{c}\right)\lesssim\Delta^{1-\beta/2}.
\]
It remains to bound 
$\mathbb{P}\left[\left|\Delta T_{k\Delta}^{(1)}\right|\geq2\xi_{0}\sqrt{\Delta}\ln(n)\right]$.
We have that: 
\[
\mathbb{P}\left[\left|\Delta T_{k\Delta}^{(1)}\right|\geq2\xi_{0}\sqrt{\Delta}\ln(n)\right]
\leq2\mathbb{P}\left[\exp\left(\varepsilon\int_{k\Delta}^{(k+1)\Delta}\xi(X_{s^{-}})dL_{s}^{(1)}\right)\geq n^{2\ep\xi_{0}\sqrt{\Delta}}\right].
\]
By Proposition \ref{prop_martingale_exponentielle}, for any $\varepsilon$,
\[
D_{t}:=\exp\left(\varepsilon\int_{k\Delta}^{t}\xi(X_{s^{-}})dL_{s}^{(1)}-\int_{k\Delta}^{t}\int_{\left|z\right|
\leq\sqrt{\Delta}}\left(\exp(\varepsilon z\xi(X_{s^{-}})-1-\varepsilon z\xi(X_{s^{-}})\right)\nu(dz)\right)
\]
 is a local martingale. Let us set $\varepsilon=1/(2\xi_{0}\Delta^{1/2}$).
There exists an increasing sequence of stopping times $\tau_{N}$
such that, for any $N$, 
\begin{eqnarray*}
F & := & \mathbb{P}\left[\exp\left(\frac{1}{2\xi_{0}\Delta^{1/2}}\int_{k\Delta}^{(k+1)\Delta\wedge\tau_{N}}\xi(X_{s^{-}})dL_{s}^{(1)}\right)\geq n\right]\\
 & \leq & 
n^{-1}\E\left(\exp\left(\int_{k\Delta}^{(k+1)\Delta\wedge\tau_{N}}\int_{\left|z\right|\leq\sqrt{\Delta}}
\left(\exp\left(\frac{z\xi(X_{s^{-}})}{2\xi_{0}\Delta^{1/2}}\right)-1-\frac{z\xi(X_{s^{-}})}{2\xi_{0}\Delta^{1/2}}\right)\nu(dz)\right)\right)\\
 & \leq & 
n^{-1}\exp\left(2\Delta\int_{\vert z\vert\leq\sqrt{\Delta}}\frac{\xi_{0}^{2}z^{2}}{4\xi_{0}^{2}\Delta}\nu(dz)\right)
\leq n^{-1}\exp\left(\int_{\mathbb{R}}z^{2}\nu(dz)\right)\leq n^{-1}.
\end{eqnarray*}
When $N\rightarrow\infty$, by dominated convergence, we obtain: 
\begin{equation}
\mathbb{P}\left(\left|\Delta T_{k\Delta}^{(1)}\right|\geq\xi_{0}\sqrt{\Delta}\ln(n)\right)\lesssim n^{-1}.\label{eq:majoration_T1}
\end{equation}

\paragraph{Bound of $\mathbb{P}\left(\Omega_{X,k}\cap\Omega_{N,k}^{c}\cap(X_{k\Delta}\in A)\right)$. }

We recall that \\
 $N'_{k\Delta}=\mu\left(]k\Delta,(k+1)\Delta],\:[-\Delta^{1/4},\Delta^{1/4}]^{c}\right)$.
We have: 

\[
\Omega_{N,k}^{c}=\left\{ N_{k\Delta}^{'}=1\right\} \cup\left\{ N_{k\Delta}^{'}\geq2\right\} 
\]
with 
\[
\mathbb{P}\left(N_{k\Delta}^{'}=1\right)\lesssim\Delta^{1-\beta/4}
\quad\textrm{and}\quad
\mathbb{P}\left(N_{k\Delta}^{'}\geq2\right)\lesssim\Delta^{2-\beta/2}.
\]
Then
 $\mathbb{P}\left(\Omega_{N,k}^{c}\cap\left\{ N_{k\Delta}^{'}\geq2\right\} \right)\lesssim\Delta^{2-\beta/2}.$
We can write: 
\begin{eqnarray*}
G & := & \mathbb{P}\left(\Omega_{X,k}\cap(X_{k\Delta}\in A)\cap(N_{k\Delta}^{'}=1)\right)\\
 & \leq & \mathbb{P}\left(N_{k\Delta}^{'}=1\right)\mathbb{P}\left(\left.\left|\Delta T_{k\Delta}^{(2)}+\Delta T_{k\Delta}^{(3)}\right|
\leq2C_{\Delta}\right|N'_{k\Delta}=1\right)\\
 & + & \mathbb{P}\left(N_{k\Delta}^{'}=1\right)\mathbb{P}\left(\left\{ \left.\left|\Delta T_{k\Delta}^{(2)}+\Delta T_{k\Delta}^{(3)}\right|
\geq2C_{\Delta}\right|N'_{k\Delta}=1\right\} \:\cap\Omega_{X,k}\cap(X_{k\Delta}\in A)\right).
\end{eqnarray*}
By \eqref{eq:majoration_I}, \eqref{eq:majoration_Z} and \eqref{eq:majoration_T1},
we obtain: 
\begin{eqnarray*}
H & := & \mathbb{P}\left(\left\{ \left.\left|\Delta T_{k\Delta}^{(2)}+\Delta T_{k\Delta}^{(3)}\right|\geq2C_{\Delta}\right|N'_{k\Delta}=1\right\} 
\:\cap\Omega_{X,k}\cap(X_{k\Delta}\in A)\right)\\
 & \leq & \mathbb{P}\left(\Delta\left|b_{A}(X_{k\Delta})+I_{k\Delta}+Z_{k\Delta}+T_{k\Delta}^{(1)}\right|>C_{\Delta}\right)\\
 & \lesssim & \Delta+n^{-1}.
\end{eqnarray*}
It remains to bound 
$J:=\mathbb{P}\left(\left|\Delta T_{k\Delta}^{(2)}+\Delta T_{k\Delta}^{(3)}\right|\leq\left.2C_{\Delta}\right|N_{k\Delta}^{'}=1\right)$.
If 
$N_{k\Delta}^{'}=1$, then $\left|\Delta T_{k\Delta}^{(3)}\right|=\vert\int_{k\Delta}^{(k+1)\Delta}\xi(X_{s^{-}})dL_{s}^{(3)}\vert\geq\xi_{1}\Delta^{1/4}$.
Then 
$J\leq\mathbb{P}\left(\Delta\left|T_{k\Delta}^{(2)}\right|\geq\xi_{1}\Delta^{1/4}-2C_{\Delta}\right)$.
Let us set 
$n_{0}=\left\lceil \frac{1}{1-\beta/2}\right\rceil $ 
and
$a=\left(\xi_{0}n_{0}\right)^{-1}\left(\xi_{1}\Delta^{1/4}-2C_{\Delta}\right).$
We have: 
\begin{eqnarray*}
J & \leq & \mathbb{P}\left[\mu(]k\Delta,(k+1)\Delta],[-a,a]^{c})\geq1\right]\\
&+&\mathbb{P}\left[\mu(]k\Delta,(k+1)\Delta],\:[-a,-\Delta^{1/2}]\cup[\Delta^{1/2},a])\geq n_{0}\right]\\
 & \leq & \Delta\nu([-a,a]^{c})+\Delta^{n_{0}}\nu([-\Delta^{1/2},\Delta^{1/2}]^{c})^{n_{0}}\\
 & \lesssim & \Delta^{1-\beta/4}+\Delta.
\end{eqnarray*}
Then $\mathbb{P}(\Omega_{X,k}\cap\Omega_{N,k}^{c})\leq\mathbb{P}(N'_{k\Delta}=1)\Delta^{1-\beta/4}+\mathbb{P}(N'_{k\Delta}=2)\lesssim\Delta^{2-\beta/2}.$

\paragraph{Bound of $\left(\E\left[\left.\left(Z_{k\Delta}+T_{k\Delta}\right)\units{\Omega_{X,k}\cap\Omega_{N,k}\cap(X_{k\Delta}\in A)}\right|
\rond{F}_{k\Delta}\right]\right)^{2}$. }

\subparagraph{If $\sigma$ and $\xi$ are constants. }

Let us set 
$E:=\left(\E\left[\left.\left(Z_{k\Delta}+T_{k\Delta}\right)
\units{\Omega_{X,k}\cap\Omega_{N,k}\cap(X_{k\Delta}\in A)}\right|\rond{F}_{k\Delta}\right]\right)^{2}$
and 
\begin{eqnarray*}
\Omega_{I,k} & = & \left\{ \omega,\left|I_{k\Delta}\right|\leq1,\cap\left|\Delta Z_{k\Delta}\right|
\leq\sigma_{0}\sqrt{\Delta}\ln(n),\cap\left|\Delta T_{k\Delta}^{(1)}\right|\leq2\xi_{0}\sqrt{\Delta}\ln(n)\right\} .
\end{eqnarray*}
 By \eqref{eq:majoration_I}, \eqref{eq:majoration_Z} and \eqref{eq:majoration_T1},
$\mathbb{P}\left(\Omega_{I,k}^{c}\right)\leq\Delta+n^{-1}.$ 
Then,
by a Markov inequality:
 \[
E\lesssim\Delta\ln^{2}(n)+\left(\E\left[\left.\left(Z_{k\Delta}+T_{k\Delta}\right)
\units{\Omega_{X,k}\cap\Omega_{N,k}\cap\Omega_{I,k}\cap(X_{k\Delta}\in A)}\right|\rond{F}_{k\Delta}\right]\right)^{2}.
\]
Let us introduce the set 
$\Omega_{ZT,k}:=\left\{ \omega,\left|Z_{k\Delta}+T_{k\Delta}\right|\leq C_{\Delta}\Delta^{-1}-b_{max}-1\right\} $.
On 
$\Omega_{I,k}$, $\left|I_{k\Delta}\right|\leq1$
and therefore: 
\[
\Omega_{ZT,k}\cap\Omega_{I,k}\subseteq\Omega_{X,k}\cap\Omega_{I,k}\subseteq
\left\{ \omega,\left|Z_{k\Delta}+T_{k\Delta}\right|\leq C_{\Delta}\Delta^{-1}+b_{max}+1\right\} \cap\Omega_{I,k}.
\]
Then 
\[
E\lesssim\Delta\ln^{2}(n)+F^{2}+G^{2}
\]
where 
$F=\E\left[\left.\left(Z_{k\Delta}+T_{k\Delta}\right)
\units{\Omega_{ZT,k}\cap\Omega_{N,k}\cap\Omega_{I,k}\cap(X_{k\Delta}\in A)}\right|\rond{F}_{k\Delta}\right]$
and \\
$G=\E\left[\left.\left(Z_{k\Delta}+T_{k\Delta}\right)\units{\Omega_{ZT,k}^{c}\cap\Omega_{X,k}\cap\Omega_{N,k}\cap\Omega_{I,k}\cap(X_{k\Delta}\in A)}\right|
\rond{F}_{k\Delta}\right]$.
As $\sigma$ and $\xi$ are constants, the terms 
\[
Z_{k\Delta}=\frac{\sigma_{0}}{\Delta}\int_{k\Delta}^{(k+1)\Delta}dW_{s}
\quad\textrm{and}\quad 
T_{k\Delta}=\frac{\xi_{0}}{\Delta}\int_{k\Delta}^{(k+1)\Delta}dL_{s}
\]
are centred and independent. Then $F=0$. Moreover, on $\Omega_{N,k}$,
$T_{k\Delta}^{(3)}=0$. Then 
\[
\left|G\right|\lesssim\left|\E\left[\left.\left(Z_{k\Delta}+T_{k\Delta}^{(1)}+T_{k\Delta}^{(2)}\right)
\units{\Omega_{X,k}\cap\Omega_{ZT,k}^{c}\cap\Omega_{N,k}\cap\Omega_{I,k}\cap(X_{k\Delta}\in A)}\right|\rond{F}_{k\Delta}\right]\right|.
\]
Let us set $c_{b}=b_{max}+1$. On $\Omega_{I,k}\cap\Omega_{X,k}$, $\left|Z_{k\Delta}+T_{k\Delta}^{(1)}+T_{k\Delta}^{(2)}\right|\lesssim\ln(n)\Delta^{-1/2}$,
and 
\begin{eqnarray*}
\left|G\right| & \lesssim & 
\frac{\ln(n)}{\sqrt{\Delta}}\left(\mathbb{P}\left(\left|Z_{k\Delta}+T_{k\Delta}^{(1)}+T_{k\Delta}^{(2)}\right|
\in\left[C_{\Delta}\Delta^{-1}-c_b,C_{\Delta}\Delta^{-1}+c_b\right]\units{\Omega_{I,k}}\right)\right)\\
 & = & 2\frac{\ln(n)}{\sqrt{\Delta}}\int_{\mathbb{R}}\mathbb{P}\left(T_{k\Delta}^{(2)}
\in\left[C_{\Delta}\Delta^{-1}-c_b-x,C_{\Delta}\Delta^{-1}+c_b-x\right]\units{\Omega_{I,k}}\right)\\
 & \times & \mathbb{P}\left(\left.Z_{k\Delta}+T_{k\Delta}^{(1)}\in dx\right|T_{k\Delta}^{(2)}
\in\left[C_{\Delta}\Delta^{-1}-c_b-x,C_{\Delta}\Delta^{-1}+c_b-x\right]\units{\Omega_{I,k}}\right).
\end{eqnarray*}
On $\Omega_{I,k}$, $\left|Z_{k\Delta}+T_{k\Delta}^{(1)}\right|\leq(\sigma_{0}+2\xi_{0})\ln(n)\Delta^{-1/2}$.
 Then 
\begin{equation}
\left|G\right|  \lesssim  
\frac{\ln(n)}{\sqrt{\Delta}}\left[\sup_{C\geq\xi_{0}\ln(n)\Delta^{-1/2}}\mathbb{P}\left(T_{k\Delta}^{(2)}
\in\left[C,C+2c_b\right]\right)\right].
\label{eq:majoration_E}
\end{equation}
We recall that $L_{t}^{(2)}$ is a compound Poisson process in which
all the jumps are greater than $\sqrt{\Delta}$ and smaller than $\Delta^{1/4}$.
Let us denote by $\tau_{i}$ the times of the jumps of size in $[\sqrt{\Delta},\Delta^{1/4}]$
and by $\zeta_{i}$ the size of the jumps. We set 
$a_{j}=\xi_{0}^{-1}C\Delta-\sum_{i=1}^{j-1}\zeta_{i}$
and 
$c:=\xi_{0}^{-1}(2b_{max}+2)$. 
Then, as $\xi$ is constant equal
to $\xi_{0}$: 
\begin{eqnarray*}
H & := & \mathbb{P}\left(T_{k\Delta}^{(2)}\in\left[C,C+2b_{max}+2\right]\right)\\
 & \leq & \sum_{j=1}^{\infty}\mathbb{P}\left(j\textrm{ jumps }\geq\sqrt{\Delta},\:\textrm{last jump }\in\left[a_{j},a_{j}+c\Delta\right]\right)\\
 & \lesssim & 2\sup_{a\geq\sqrt{\Delta}}\mathbb{P}\left(\textrm{1 jump }\in\left[a,a+c\Delta\right]\right)
=2\Delta\sup_{a\geq\sqrt{\Delta}}\nu\left(\left[a,a+c\Delta\right]\right).
\end{eqnarray*}
By A\ref{hypo_nu_reguliere}, 
\begin{equation}
H\lesssim\Delta\sup_{a\geq\sqrt{\Delta}}\left[\frac{1}{a^{\beta}}-\frac{1}{\left(a+c\Delta\right)^{\beta}}\right]
\lesssim\sqrt{\Delta}\Delta^{1-\beta/2}\label{eq:majoration_F}
\end{equation}
and, by \eqref{eq:majoration_E} and \eqref{eq:majoration_F}, 
\[
E\lesssim\Delta\ln^{2}(n)+\frac{\ln^{2}(n)}{\Delta}\Delta\Delta^{2-\beta}\lesssim\Delta\ln^{2}(n)+\Delta^{2-\beta}\ln^{2}(n).
\]

\begin{req}\label{ref:nu_pas_continue}

If $\nu$ is not absolutely continuous, inequality \ref{eq:majoration_F} is not valid.  
We obtain: 
\[ H\lesssim 2\Delta \sup_{a\geq \sqrt{\Delta}}\nu([a,a+c\Delta])\lesssim \Delta^{1-\beta/2}\]
Therefore
\[ 
E\leq\Delta\ln^{2}(n)+G^{2}\lesssim\Delta\ln^{2}(n)+\Delta^{1-\beta}\ln^{2}(n).\]

\end{req}

\subparagraph{If $\sigma$ or $\xi$ are not constants. }

The problem is that $Z_{k\Delta}$ and $T_{k\Delta}$ are not symmetric
and we can't apply directly the previous method. We replace them by
two centred terms. The following lemma is very useful. 

\begin{lem}\label{lem_psi}

Let $f$ be a $\rond{C}^{2}$ function such that $f$ and $f'$ are
Lipschitz. Let us set, for any $t\in]k\Delta,(k+1)\Delta]$: 
\[
\psi_{f}(X_{k\Delta},t)=f'(X_{k\Delta})\left(\sigma(X_{k\Delta})\int_{k\Delta}^{t}dW_{s}+\xi(X_{k\Delta})\int_{k\Delta}^{t}z\tilde{\mu}(ds,dz)\right).
\]
We have:  
\[
\E\left[\left(f(X_{t})-f(X_{k\Delta})-\psi_{f}(X_{k\Delta},t)\right)^{2}\units{\Omega_{N,k}}\units{X_{k\Delta}\in A}\right]\lesssim\Delta^{2-\beta/4}.
\]

\end{lem}

Lemma 4 is proved below. Let us set 
\[
\bar{Z}_{k\Delta}=\frac{1}{\Delta}\int_{k\Delta}^{(k+1)\Delta}\left(\sigma(X_{k\Delta})+\psi_{\sigma}(X_{k\Delta,s})\right)dW_{s},
\]
 \[
\bar{T}_{k\Delta}^{(i)}=\frac{1}{\Delta}\int_{k\Delta}^{(k+1)\Delta}\left(\xi(X_{k\Delta})+\psi_{\xi}(X_{k\Delta,s})\right)dL_{s}^{(i)}
\quad\textrm{and}\quad
\bar{T}_{k\Delta}=\bar{T}_{k\Delta}^{(1)}+\bar{T}_{k\Delta}^{(2)}+\bar{T}_{k\Delta}^{(3)}.
\]
The terms $\bar{Z}_{k\Delta}$ and $\bar{T}_{k\Delta}$ are symmetric.
By lemma \ref{lem_psi}, 
\begin{eqnarray}
\E\left[\left(\bar{Z}_{k\Delta}-Z_{k\Delta}\right)^{2}\units{\Omega_{N,k}}\units{X_{k\Delta}\in A}\right] 
& = & 
\frac{1}{\Delta^{2}}\E\left[\int_{k\Delta}^{(k+1)\Delta}\left(\sigma(X_{s})-\sigma(X_{k\Delta})-\psi_{\sigma}(X_{k\Delta,s})\right)^{2}ds\right]\nonumber \\
 & \lesssim & \Delta^{1-\beta/4}.
\label{eq:majoration_Zbar_Z}
\end{eqnarray}
We prove in the same way that 
\begin{equation}
\E\left[\left(\bar{T}_{k\Delta}-T_{k\Delta}\right)^{2}\units{\Omega_{N,k}}\units{X_{k\Delta}\in A}\right]\leq\Delta^{1-\beta/4}.
\label{eq:majoration_Tbar_T}
\end{equation}
Let us set $U_{k\Delta}=\Delta^{-1}\xi(X_{k\Delta^{-}})\int_{k\Delta}^{(k+1)\Delta}dL_{s}^{(2)}$.
By Result \ref{res_burkholder} and Proposition \ref{prop_gloter},
\begin{equation}
\E\left[\Delta^{2}\left(\bar{T}_{k\Delta}^{(2)}-U_{k\Delta}\right)^{2}\right]
=\E\left[\int_{k\Delta}^{(k+1)\Delta}\int_{\mathbb{R}}\left(\psi_{\xi}(X_{k\Delta,s})\right)^{2}z^{2}\nu(dz)ds\right]\leq\Delta^{2-\beta/4}.
\label{eq:majoration_U_Tbar}
\end{equation}
 Let us introduce the set 
\begin{eqnarray*}
\bar{\Omega}_{I,k} & = &
 \left\{ \omega,\left|I_{k\Delta}\right|+\left|Z_{k\Delta}-\bar{Z}_{k\Delta}\right|+\left|T_{k\Delta}-\bar{T}_{k\Delta}\right|\leq3\right\} \\
 & \bigcap & \left\{ \left|\Delta\bar{Z}_{k\Delta}\right|\leq\sigma_{0}\sqrt{\Delta}\ln(n)+\Delta,\:\left|\Delta\bar{T}_{k\Delta}^{(1)}\right|\leq2\xi_{0}\sqrt{\Delta}\ln(n)+\Delta\right\} \\
 & \bigcap & \left\{ \left|\Delta(\bar{T}_{k\Delta}^{(2)}-U_{k\Delta})\right|\leq\xi_{0}\sqrt{\Delta}\right\} .
\end{eqnarray*}
By \eqref{eq:majoration_I}, \eqref{eq:majoration_Z}, \eqref{eq:majoration_T1},
\eqref{eq:majoration_Zbar_Z}, \eqref{eq:majoration_Tbar_T}, \eqref{eq:majoration_U_Tbar}
and Markov inequalities, we obtain:

\begin{equation}
\mathbb{P}\left(\bar{\Omega}_{I,k}^{c}\right)\lesssim\Delta^{1-\beta/4}+\frac{1}{n}.
\label{eq:majoration_omegaI}
\end{equation}
Then
\begin{eqnarray}
E & := & 
\left(\E\left[\left.\left(Z_{k\Delta}+T_{k\Delta}\right)\units{\Omega_{X,k}\cap\Omega_{N,k}\cap(X_{k\Delta}\in A)}\right|
\rond{F}_{k\Delta}\right]\right)^{2}
\label{eq:def_E}\\
 & \lesssim & \Delta^{1-\beta/2}\ln^{2}(n)+\left(\E\left[\left.\left(\bar{Z}_{k\Delta}+\bar{T}_{k\Delta}\right)
\units{\Omega_{X,k}\cap\Omega_{N,k}\cap(X_{k\Delta}\in A)\cap\bar{\Omega}_{I,k}}\right|\rond{F}_{k\Delta}\right]\right)^{2}.\nonumber 
\end{eqnarray}
Let us introduce the set:
 \[
\bar{\Omega}_{ZT,k}:=\left\{ \omega,\left|\bar{Z}_{k\Delta}+\bar{T}_{k\Delta}\right|\leq C_{\Delta}\Delta^{-1}-b_{max}-3\right\} .
\]
We have that 
\[
\bar{\Omega}_{ZT,k}\cap\bar{\Omega}_{I,k}\subseteq\Omega_{X,k}\cap\bar{\Omega}_{I,k}
\subseteq\left\{ \omega,\left|\bar{Z}_{k\Delta}+\bar{T}_{k\Delta}\right|\leq C_{\Delta}\Delta^{-1}+b_{max}+3\right\} \cap\bar{\Omega}_{I,k}.
\]
Given the filtration $\rond{F}_{k\Delta}$, the sum $\bar{Z}_{k\Delta}+\bar{T}_{k\Delta}$
is symmetric. Then 
\[
\E\left[\left.\left(\bar{Z}_{k\Delta}+\bar{T}_{k\Delta}\right)
\units{\bar{\Omega}_{ZT,k}\cap\Omega_{N,k}\cap(X_{k\Delta}\in A)}\right|\rond{F}_{k\Delta}\right]=0.
\]
 Moreover, on $\Omega_{N,k}$, $\bar{T}_{k\Delta}^{(3)}=0$. Then,
by \eqref{eq:def_E}, 
\[
E\lesssim\Delta^{1-\beta/2}\ln^{2}(n)+G^{2}+H^{2}
\]
where 
$G:=\E\left[\left.\left(\bar{Z}_{k\Delta}+\bar{T}_{k\Delta}^{(1)}+\bar{T}_{k\Delta}^{(2)}\right)
\units{\Omega_{X,k}\cap\Omega_{ZT,k}^{c}\cap\Omega_{N,k}\cap\Omega_{I,k}\cap(X_{k\Delta}\in A)}\right|\rond{F}_{k\Delta}\right]$
and 
$H:=\E\left[\left.\left(\bar{Z}_{k\Delta}+\bar{T}_{k\Delta}^{(1)}+\bar{T}_{k\Delta}^{(2)}\right)
\units{\Omega_{X,k}\cap\Omega_{ZT,k}\cap\Omega_{N,k}\cap\Omega_{I,k}^{c}\cap(X_{k\Delta}\in A)}\right|\rond{F}_{k\Delta}\right]$.
We have that 
$H^{2}\lesssim\Delta^{-1}\ln^{2}(n)\mathbb{P}^{2}(\Omega_{I,k}^{c})\lesssim\Delta^{1-\beta/2}\ln^{2}(n)$.
The end of the proof is the same as in the case of $\sigma$ and $\xi$
constants. We obtain that 
\[
\left| G\right|  \lesssim 
\frac{\ln(n)}{\sqrt{\Delta}}\sup_{C\geq\kappa_{0}\ln(n)\Delta^{-1/2}}\mathbb{P}\left(U_{k\Delta}\in\left[C,C+2b_{max}+6\right]\right)
\lesssim\sqrt{\Delta}\Delta^{1-\beta/2}.
\]

\subsubsection{Proof of Lemma \ref{lem_psi}}

According to the Itô formula (see for instance \citet{applebaum},
Theorem 4.4.7 p251), we have that 
\[
f(X_{t})-f(X_{k\Delta})=I_{1}+I_{2}+I_{3}+I_{4}
\]
where \begin{gather*}
I_{1}=\int_{k\Delta}^{t}f'(X_{s})\sigma(X_{s})dW_{s},\quad I_{2}=
\int_{k\Delta}^{t}\int_{\mathbb{R}}\left(f\left(X_{s^{-}}+z\xi(X_{s^{-}})\right)-f(X_{s^{-}})\right)\tilde{\mu}(ds,dz)\\
I_{3}=\int_{k\Delta}^{t}\int_{z\in\mathbb{R}}\left[f(X_{s}+z\xi(X_{s}))-f(X_{s})-z\xi(X_{s})f'(X_{s})\right]\nu(dz)ds\\
I_{4}=\int_{k\Delta}^{t}\left[f'(X_{s})b(X_{s})+f''(X_{s})\sigma^{2}(X_{s})/2\right]ds.
\end{gather*}
By Proposition \ref{prop_gloter}, for any $t\leq(k+1)\Delta$, we
have: 
\begin{eqnarray*}
Q&:=&\E\left[\left(I_{1}-f'(X_{k\Delta})\sigma(X_{k\Delta})\int_{k\Delta}^{t}dW_{s}\right)^{2}\right] \\
& = & 
\E\left[\left(\int_{k\Delta}^{t}\left(\sigma(X_{s})f'(X_{s})-\sigma(X_{k\Delta})f'(X_{k\Delta})\right)dW_{s}\right)^{2}\right]\\
 & = & \int_{k\Delta}^{t}\left(\sigma(X_{s})f'(X_{s})-\sigma(X_{k\Delta})f'(X_{k\Delta})\right)^{2}ds\lesssim\Delta^{2}.
\end{eqnarray*}
We can write: 
\begin{eqnarray*}
E & := & \E\left[\left(I_{2}-f'(X_{k\Delta})\xi(X_{k\Delta^{-}})\int_{k\Delta}^{t}dL_{s}^{(1)}+dL_{s}^{(2)}\right)^{2}\units{\Omega_{N,k}}\right]\\
 & \leq & 2\int_{k\Delta}^{t}\int_{\vert z\vert\leq\Delta^{1/4}}\E\left[\left(f(X_{s}+z\xi(X_{s}))-f(X_{s})-z\xi(X_{s})f'(X_{s})\right)^{2}\right]\nu(dz)ds\\
 & + & 2\int_{k\Delta}^{t}\int_{\vert z\vert\leq\Delta^{1/4}}\E\left[z^{2}\left(\xi(X_{s})f'(X_{s})-\xi(X_{k\Delta})f'(X_{k\Delta})\right)^{2}\right]\nu(dz)ds.
\end{eqnarray*}
The function $f$ is $\rond{C}^{2}$, then, by the Taylor formula,
for any $s\in[k\Delta,t]$, $z\in\mathbb{R}$, there exists $\zeta_{s,z}$
in $[X_{s},X_{s}+z\xi(X_{s})]$ such that: 
\[
f\left(X_{s}+z\xi(X_{s})\right)-f(X_{s})-z\xi(X_{s})f'(X_{s})=\frac{z^{2}\xi^{2}(X_{s})}{2}f''(\zeta_{s,z}).
\]
Then, as $\xi$ and $f''$ are bounded: 
\[
\E\left[\left(f(X_{s}+z\xi(X_{s})-f(X_{s})-z\xi(X_{s})f'(X_{s})\right)^{2}\right]=
\frac{z^{4}}{4}\E\left[\left(\xi(X_{s})f''(\zeta_{s,z})\right)^{2}\right]\lesssim z^{4}
\]
and, by Result \ref{result_blumenthal}, for any $t\leq(k+1)\Delta$,
\begin{eqnarray*}
F & := & \int_{k\Delta}^{t}\int_{\vert z\vert\leq\Delta^{1/4}}\E\left[\left(f(X_{s}+z\xi(X_{s})-f(X_{s})-z\xi(X_{s})f'(X_{s})\right)^{2}\right]\nu(dz)ds\\
 & \lesssim & \Delta\int_{\vert z\vert\leq\Delta^{1/4}}z^{4}\nu(dz)\lesssim\Delta^{2-\beta/4}.
\end{eqnarray*}
The functions $\xi$ and $f'$ are Lipschitz, then by Proposition
\ref{prop_gloter}, 
\[
\E\left[z^{2}\left(\xi(X_{s})f'(X_{s})-\xi(X_{k\Delta})f'(X_{k\Delta})\right)^{2}\right]
\lesssim z^{2}\E\left[\left(X_{s}-X_{k\Delta}\right)^{2}\right]
\lesssim\Delta z^{2}
\]
 and consequently, for any $t\leq(k+1)\Delta$:
\[
\int_{k\Delta}^{t}\int_{\vert z\vert\leq\Delta^{1/4}}\E\left[z^{2}\left(\xi(X_{s})f'(X_{s})-\xi(X_{k\Delta})f'(X_{k\Delta})\right)^{2}\right]\nu(dz)ds
\lesssim\Delta^{2-\beta/4}
\]
then $E\lesssim\Delta^{2-\beta/4}$. By the same way, we obtain that
\[
\E\left[I_{3}^{2}\right]\leq
\E\left[\int_{k\Delta}^{t}\int_{\vert z\vert\leq\Delta^{1/4}}\left(\frac{z^{2}\xi^{2}(X_{s})}{2}f''(\zeta_{s,z})\right)^{2}\nu(dz)ds\right]
\lesssim \Delta^{2-\beta/4}.
\]
The functions $b$ and $f'$ are Lipschitz and $f''$ and $\sigma$
are bounded, then, for any $t\leq(k+1)\Delta$ : 
\[
\E\left[I_{4}^{2}\right]\lesssim\Delta\int_{k\Delta}^{t}\left(1+\E\left[X_{s}^{4}\right]\right)ds\lesssim\Delta^{2}.
\]
Then, for any $t\leq(k+1)\Delta$: 
\[
\E\left[\left(f(X_{t})-f(X_{k\Delta})-\psi_{f}(X_{k\Delta},t)\right)\right]\leq\Delta^{2-\beta/4}.
\]

\subsection{Proof of Theorem \ref{thme_sauts_coupes_adaptif}}

As previously, we only bound the risk on $\Omega_{n}$. As in Subsection 
\ref{sub:Proof-of-Theorem_adaptatif}, we introduce the function $p(m,m')$
such that $p(m,m')=12(pen(m)+pen(m'))$. On $\Omega_{n}$, for any
$m\in\rond{M}_{n}$, we have: 
\begin{eqnarray*}
\left\Vert \tilde{b}_{\tilde{m}}-b_{A}\right\Vert _{n}^{2} & \leq & 
3\left\Vert b_{m}-b_{A}\right\Vert _{n}^{2}+\frac{224}{n}\sum_{k=1}^{n}b_{A}^{2}(X_{k\Delta})\units{\Omega_{X,k}^{c}}+
I_{k\Delta}^{2}+2\left(Z_{k\Delta}^{2}+T_{k\Delta}^{2}\right)\units{\Omega_{X,k}\cap\Omega_{Z,k}^{c}}\\
 & + & \frac{224}{n}\sum_{k=1}^{n}\left(\E\left[\left.\left(Z_{k\Delta}+
T_{k\Delta}\right)\units{\Omega_{X,k}\cap\Omega_{Z,k}}\right|\rond{F}_{k\Delta}\right]\right)^{2}\\
 & + & 24\sup_{t\in\rond{B}_{m,\hat{m}}}\left(\tilde{\nu}_{n}^{2}(t)-p(m,\tilde{m})\right)+4pen(m).
\end{eqnarray*}
 It remains only to bound
\[
\E\left[\sup_{t\in\rond{B}_{m,\hat{m}}}\left(\tilde{\nu}_{n}^{2}(t)-p(m,\tilde{m})\right)\right]
\leq\sum_{m'}\E\left[\sup_{t\in\rond{B}_{m,m'}}\left(\tilde{\nu}_{n}^{2}(t)-p(m,\tilde{m})\right)\right].
\]

As in the proof of Theorem \ref{thme_adaptatif}, we bound the quantity
\[
\E\left[\left.\exp\left(\varepsilon t(X_{k\Delta})\left(\tilde{Z}_{k\Delta}+\tilde{T}_{k\Delta}\right)\right)\right|\rond{F}_{k\Delta}\right].
\]
We have that
 \[
\E\left[\left.\exp\left(\varepsilon t(X_{k\Delta})Z_{k\Delta}\right)\units{\Omega_{N,k}}\right|\rond{F}_{k\Delta}\right]\leq
\exp\left(\frac{\varepsilon^{2}\sigma_{0}^{2}t^{2}(X_{k\Delta})}{2\Delta}\right).
\]
The truncated Lévy process 
$\tilde{L}_{t}=\int_{0}^{t}\int_{\left|z\right|\leq\Delta^{1/4}}z\tilde{\mu}(ds,dz)$
satisfies Assumption A\ref{hypo_loi_xi} and then there exists a constant
$c$ such that: 
\[
\E\left[\left.\exp\left(\varepsilon t(X_{k\Delta})T_{k\Delta}\right)\units{\Omega_{N,k}}\right|\rond{F}_{k\Delta}\right]
\leq\exp\left(\frac{c\varepsilon^{2}\xi_{0}^{2}t^{2}(X_{k\Delta})}{\Delta\left(1-\varepsilon/\varepsilon_{1}\right)}\right).
\]
As $Z_{k\Delta}\units{\Omega_{N,k}}$ and $T_{k\Delta}\units{\Omega_{N,k}}$
are centred, we obtain: 
\[
\E\left[\left.\exp\left(\varepsilon\left|t(X_{k\Delta})\left(Z_{k\Delta}+T_{k\Delta}\right)\right|\right)\units{\Omega_{N,k}}\right|\rond{F}_{k\Delta}\right]
\leq 2\exp\left(\frac{c\varepsilon^{2}\left(\sigma_{0}^{2}+\xi_{0}^{2}\right)t^{2}(X_{k\Delta})}{\Delta\left(1-\varepsilon/\varepsilon_{1}\right)}\right)
\]
and then 
\[
\E\left[\left.\exp\left(\varepsilon\left|t(X_{k\Delta})\left(\tilde{Z}_{k\Delta}+\tilde{T}_{k\Delta}\right)\right|\right)
\units{\Omega_{N,k}\cap\Omega_{X,k}}\right|\rond{F}_{k\Delta}\right]
\leq2\exp\left(\frac{c\varepsilon^{2}\left(\sigma_{0}^{2}+\xi_{0}^{2}\right)t^{2}(X_{k\Delta})}{\Delta\left(1-\varepsilon/\varepsilon_{1}\right)}\right).
\]
We conclude as in the proof of Theorem \ref{thme_adaptatif}. 

\begin{table}[p]
\caption{Model 1: Ornstein-Uhlenbeck and binomial law}
\label{Flo:table_Model1}

\[
b(x)=-2x,\;\sigma(x)=\xi(x)=1\:\textrm{and compound Poisson process (binomial law)}
\]

\begin{tabular}{|c|c||c|c|c|c||c|c|c|c|}
\hline 
 &  & \multicolumn{4}{c||}{first estimator} & \multicolumn{4}{c|}{truncated estimator }\tabularnewline
\hline 
$n$ & $\Delta$ & $\hat{m}_{a}$ & $\hat{r}_{a}$ & $risk_{1}$ & $or_{1}$ & $\tilde{m}_{a}$ & $\tilde{r}_{a}$ & $risk_{2}$ & $or_{2}$\tabularnewline
\hline
\hline 
$10^{3}$ & $10^{-1}$ & 0 & 1.02 & 0.044 & 1.3 & 0 & 1.02 & 0.044 & 1.3\tabularnewline
\hline 
$10^{4}$ & $10^{-1}$ & 0 & 1.02 & 0.011 & 1.3 & 0 & 1.02 & 0.011 & 1.3\tabularnewline
\hline 
$10^{3}$ & $10^{-2}$ & 0 & 1.02 & 0.55 & 1.04 & 0 & 1.02 & 0.55 & 1.04\tabularnewline
\hline 
$10^{4}$ & $10^{-2}$ & 0 & 1 & 0.047 & 1 & 0 & 1 & 0.047 & 1\tabularnewline
\hline 
$5.10^{4}$ & $10^{-2}$ & 0.04 & 1 & 0.010 & 1.4 & 0 & 1 & 0.0053 & 1\tabularnewline
\hline
\end{tabular}\\
$\hat{m}_{a}$, $\hat{r}_{a}$ and $\tilde{m}_{a}$, $\tilde{r}_{a}$
: average values of $\hat{m}$, $\hat{r}$ and $\tilde{m}$, $\tilde{r}$
on the 50 simulations. \\
$risk_{1}$ and $risk_{2}$ : means of the empirical errors of
the adaptive estimators. 

$or_{1}$ and $or_{2}$: means of $oracle=$empirical error of the
adaptive estimator / empirical error of the best possible estimator. 
\end{table}

\begin{table}[p]
\caption{Model 2: Double well and Laplace law}

\label{Flo:table_model2}

$b(x)=-(x-1/4)^{3}-(x+1/4)^{3}$, $\sigma(x)=\xi(x)=1$ and Laplace
law. 

\begin{tabular}{|c|c||c|c|c|c||c|c|c|c|}
\hline 
\multicolumn{1}{|c|}{} &  & \multicolumn{4}{c|}{first estimator} & \multicolumn{4}{c|}{truncated estimator }\tabularnewline
\hline 
$n$ & $\Delta$ & $\hat{m}_{a}$ & $\hat{r}_{a}$ & $risk_{1}$ & $or_{1}$ & $\tilde{m}_{a}$ & $\tilde{r}_{a}$ & $risk_{2}$ & $or_{2}$\tabularnewline
\hline
\hline 
$10^{3}$ & $10^{-1}$ & 0.02 & 1.0 & 0.12 & 3.1 & 0.02 & 1.0 & 0.12 & 3.1\tabularnewline
\hline 
$10^{4}$ & $10^{-1}$ & 1.7 & 2.1 & 2e96 & 51 & 0.4 & 2.1 & 0.04 & 1.5\tabularnewline
\hline 
$10^{3}$ & $10^{-2}$ & 0.26 & 1.2 & 1.8 & 3.1 & 0.06 & 1 & 0.51 & 1.4\tabularnewline
\hline 
$10^{4}$ & $10^{-2}$ & 0.12 & 1.5 & 0.16 & 1.8 & 0.08 & 1.2 & 0.13 & 2.4\tabularnewline
\hline 
$5.10^{4}$ & $10^{-2}$ & 0.30 & 2.5 & 0.035 & 1.6 & 0.26 & 2.5 & 0.019 & 1.8\tabularnewline
\hline
\end{tabular}

$\hat{m}_{a}$, $\hat{r}_{a}$ and $\tilde{m}_{a}$, $\tilde{r}_{a}$
: average values of $\hat{m}$, $\hat{r}$ and $\tilde{m}$, $\tilde{r}$
on the 50 simulations. \\
$risk_{1}$ and $risk_{2}$ : means of the empirical errors of
the adaptive estimators. 

$or_{1}$ and $or_{2}$: means of $oracle=$empirical error of the
adaptive estimator / empirical error of the best possible estimator. 
\end{table}

\begin{table}[p]
\caption{Model 3: Sine function and jumps not sub-exponential}
\label{Flo:table_model3}

\[
b(x)=-2x+\sin(3x),\quad\sigma(x)=\xi(x)=\sqrt{(3+x^{2})/(1+x^{2})}\quad\textrm{and}\quad\nu(dz)\propto e^{-\sqrt{az}}/\sqrt{z}dz
\]

\begin{tabular}{|c|c||c|c|c|c||c|c|c|c|}
\hline 
\multicolumn{1}{|c|}{} &  & \multicolumn{4}{c|}{first estimator} & \multicolumn{4}{c|}{truncated estimator }\tabularnewline
\hline 
$n$ & $\Delta$ & $\hat{m}_{a}$ & $\hat{r}_{a}$ & $risk_{1}$ & $or_{1}$ & $\tilde{m}_{a}$ & $\tilde{r}_{a}$ & $risk_{2}$ & $or_{2}$\tabularnewline
\hline
\hline 
$10^{3}$ & $10^{-1}$ & 0.34 & 1.2 & 0.76 & 3.6 & 0.04 & 1.2 & 0.28 & 1.9\tabularnewline
\hline 
$10^{4}$ & $10^{-1}$ & 0.8 & 2.2 & 0.082 & 1.3 & 0.68 & 2.2 & 0.073 & 1.2\tabularnewline
\hline 
$10^{3}$ & $10^{-2}$ & 0.96 & 1.2 & 18 & 6.3 & 0.02 & 1.2 & 1.3 & 1.2\tabularnewline
\hline 
$10^{4}$ & $10^{-2}$ & 0.78 & 1.4 & 1.5 & 4.3 & 0.12 & 1.4 & 0.24 & 3.3\tabularnewline
\hline 
$5.10^{4}$ & $10^{-2}$ & 0.92 & 2.3 & 0.24 & 4.3 & 0.70 & 2.3 & 0.039 & 1.3\tabularnewline
\hline
\end{tabular}

$\hat{m}_{a}$, $\hat{r}_{a}$ and $\tilde{m}_{a}$, $\tilde{r}_{a}$
: average values of $\hat{m}$, $\hat{r}$ and $\tilde{m}$, $\tilde{r}$
on the 50 simulations. \\
$risk_{1}$ and $risk_{2}$ : means of the empirical errors of
the adaptive estimators. 

$or_{1}$ and $or_{2}$: means of $oracle=$empirical error of the
adaptive estimator / empirical error of the best possible estimator. 
\end{table}

\begin{table}[p]
\caption{Model 4: Lévy process}
\label{Flo:table_model4}

\[
b(x)=-2x,\quad\sigma(x)=\xi(x)=1 \quad\textrm{and}\quad \nu(dz)= \sum_{k=0}^{\infty} 2^{k+2}(\delta_{2^{-k}}+\delta_{-2^{-k}})
\]

\begin{tabular}{|c|c||c|c|c|c||c|c|c|c|}
\hline 
\multicolumn{1}{|c|}{} &  & \multicolumn{4}{c|}{first estimator} & \multicolumn{4}{c|}{truncated estimator }\tabularnewline
\hline 
$n$ & $\Delta$ & $\hat{m}_{a}$ & $\hat{r}_{a}$ & $risk_{1}$ & $or_{1}$ & $\tilde{m}_{a}$ & $\tilde{r}_{a}$ & $risk_{2}$ & $or_{2}$\tabularnewline
\hline
\hline 
$10^{3}$ & $10^{-1}$ & 0.04& 1.06&0.110 &1.86&0.02 &1.06  &0.111  &1.95  \tabularnewline
\hline 
$10^{4}$ & $10^{-1}$ &0.06  &1.06  & 0.0172 &1.26  & 0.06 & 1.06 &0.0176 &1.22 \tabularnewline
\hline 
$10^{3}$ & $10^{-2}$ &  0.1&1.04  &1.17  &1.88  &0  &1.04  &0.61 &1.12 \tabularnewline
\hline 
$10^{4}$ & $10^{-2}$ & 0.04& 1.08 &0.11  &1.25  &0.02  &1.08  &0.068  &1.25 \tabularnewline
\hline 
$5.10^{4}$ & $10^{-2}$ & 0.08 &1.16  &0.023  &1.71  & 0 & 1.16 & 0.011 &1.09 \tabularnewline
\hline
\end{tabular}

$\hat{m}_{a}$, $\hat{r}_{a}$ and $\tilde{m}_{a}$, $\tilde{r}_{a}$
: average values of $\hat{m}$, $\hat{r}$ and $\tilde{m}$, $\tilde{r}$
on the 50 simulations. \\
$risk_{1}$ and $risk_{2}$ : means of the empirical errors of
the adaptive estimators. 

$or_{1}$ and $or_{2}$: means of $oracle=$empirical error of the
adaptive estimator / empirical error of the best possible estimator. 
\end{table}

\section{Auxiliary proofs  \label{sec:Appendice}}

\subsection{Decomposition on a lattice}

\begin{prop}\label{reseau_complique}

If there exist some constants $c_{1}$, $c_{2}$ and $K$ independent
of $D$, $n$, $\Delta$, $b$ and $\sigma$ and two constants $\alpha$
and $\beta$ independent of $n$ and $D$ such that, for any function
$t\in S_{m}+S_{m}'$:

\[
\forall\eta,\zeta>0,\;\forall t\in S_{m}+S_{m'}\;\left\Vert t\right\Vert _{\infty}\leq C\zeta,\;\Prob\left(f_{n}(t)\geq\eta,\left\Vert t\right\Vert _{n}^{2}\leq\zeta^{2}\right)\leq K\exp\left(-\frac{\eta^{2}n\beta}{\left(c_{1}\alpha^{2}\zeta^{2}+2Cc_{2}\alpha\eta\zeta\right)}\right),\]
then there exist some constants $C$ and $\kappa$ depending only
of $\nu$ such that, if $D\leq n\beta$: \[
\E\left[\sup_{t\in\rond{B}_{m,m'}}f_{n}^{2}(t)-\frac{\kappa\alpha^{2}D}{n\beta}\right]_{+}\leq CK\frac{\kappa\alpha^{2}D^{3/2}e^{-\sqrt{D}}}{n\beta}.\]

\end{prop}

Let us consider an orthonormal (for the $L_{\varpi}^{2}$-norm) basis
$\left(\psi_{\lambda}\right)_{\lambda\in\Lambda_{m,m'}}$ of $S_{m,m'}=S_{m}+S_{m'}$
such that \[
\forall\lambda,\quad\textrm{card}\left(\left\{ \lambda',\;\left\Vert \psi_{\lambda}\psi_{\lambda'}\right\Vert \neq0\right\} \right)\leq\phi_{2}.\]
Let us set \[
\bar{r}_{m,m'}=\frac{1}{\sqrt{D}}\sup_{\beta\neq0}\frac{\left\Vert \sum_{\lambda}\beta_{\lambda}\psi_{\lambda}\right\Vert _{\infty}}{\left|\beta\right|_{\infty}}.\]
We obtain that \[
\left\Vert \sum_{\lambda}\beta_{\lambda}\psi_{\lambda}\right\Vert _{\infty}\leq\phi_{2}\left|\beta\right|_{\infty}\sup_{\lambda}\left\Vert \psi_{\lambda}\right\Vert _{\infty}\quad\textrm{et}\quad\left\Vert \psi_{\lambda}\right\Vert _{\infty}\leq\sqrt{D}\left\Vert \psi_{\lambda}\right\Vert _{L^{2}}\leq\pi_{1}\sqrt{D}\left\Vert \psi_{\lambda}\right\Vert _{\varpi}\]
then \[
\bar{r}_{m,m'}\leq\bar{r}:=\phi_{2}\pi_{1}.\]
We need a lattice of which the infinite norm is bounded. We use Lemma
9 of \citet{barronbirge1999}: 

\begin{result}

There exists a $\delta_{k}$-lattice $T_{k}$ of $L_{\varpi}^{2}\cap(S_{m}+S_{m'})$
such that \[
\left|T_{k}\cap\rond{B}_{m,m'}\right|\leq\left(5/\delta^{k}\right)^{D}\]
where $\delta_{k}=2^{-k}/5$ . Let us denote by $p_{k}(u)$ the orthogonal
projection of $u$ on $T_{k}$. For any $u\in S_{m,m'}$, $\left\Vert u-p_{k}(u)\right\Vert _{\pi}\leq\delta_{k}$
and \[
\sup_{u\in p_{k}^{-1}(t)}\left\Vert u-t\right\Vert _{\infty}\leq\bar{r}_{m,m'}\delta_{k}\leq\bar{r}\delta_{k}.\]

\end{result}

Let us set $H_{k}=\ln(\left|T_{k}\cap\rond{B}_{m,m'}\right|)$. We
have that: \[
H_{k}\leq D\ln(5/\delta_{k})=D\left(k\ln(2)+\ln(5/\delta_{0})\right)\leq C(k+1)D.\]
The decomposition of $u_{k}$ on the $\delta_{k}$-lattice must be
done very carefully: the norms $\left\Vert u_{k}-u_{k-1}\right\Vert _{\varpi}$
and $\left\Vert u_{k}-u_{k-1}\right\Vert _{\infty}$ must be controlled.
Let us set \[
\rond{E}_{k}=\left\{ u_{k}\in T_{k}\cap\rond{B}_{m,m'},\quad\left\Vert u-u_{k}\right\Vert _{\varpi}\leq\delta_{k}\quad\textrm{et}\quad\left\Vert u-u_{k}\right\Vert _{\infty}\leq\bar{r}\delta_{k}\right\} .\]
We have that $\ln(\left|\rond{E}_{k}\right|)\leq H_{k}$. For any
function $u\in\rond{B}_{m,m'}$, there exist a series $(u_{k})_{k\geq0}\in\prod_{k}\rond{E}_{k}$
such that \[
u=u_{0}+\sum_{k=1}^{\infty}\left(u_{k}-u_{k-1}\right).\]
Let us consider $(\eta_{k})_{k\geq0}$ and $\eta\in\mathbb{R}$ such
that $\eta_{0}+\sum_{k=1}^{\infty}\eta_{k}\leq\eta.$ We obtain: \begin{eqnarray}
\mathbb{P}\left(\sup_{u\in\rond{B}_{m,m'}}\left|f_{n}(u)\right|>\eta\right) & \leq & \mathbb{P}\left(\exists\left(u_{k}\right)\in\prod\rond{E}_{k},\;\left|f_{n}(u_{0})+\sum_{k=1}^{\infty}f_{n}(u_{k}-u_{k-1})\right|>\eta_{0}+\sum_{k=1}^{\infty}\eta_{k}\right)\nonumber \\
 & \leq & P_{1}+\sum_{k=1}^{\infty}P_{2,k}\label{bruitsigma_eq:majorationP}\end{eqnarray}
where \[
P_{1}=\sum_{u_{0}\in\rond{E}_{0}}\mathbb{P}\left(\left|f_{n}(u_{0})\right|>\eta_{0}\right)\quad\textrm{and}\quad P_{2,k}=\sum_{u_{k}\in\rond{E}_{k}}\mathbb{P}\left(\left|f_{n}(u_{k}-u_{k-1})\right|>\eta_{k}\right).\]
As $u_{0}\in T_{0}$, $\left\Vert u_{0}\right\Vert _{\varpi}\leq1$
and $\left\Vert u_{0}\right\Vert _{\infty}\leq\bar{r}\sqrt{D}$. Moreover,
$\left\Vert u_{0}\right\Vert _{n}^{2}\leq3/2$$\left\Vert u_{0}\right\Vert _{\varpi}^{2}\leq3\delta_{0}/2$.
Then \[
\mathbb{P}\left(\left|f_{n}(u_{0})\right|>\eta_{0}\right)=\mathbb{P}\left(\left|f_{n}(u_{0})\right|>\eta_{0},\;\left\Vert u_{0}\right\Vert _{n}^{2}\leq3\delta_{0}/2\right).\]
There exist two constants $c_{1}'$ and $c_{2}'$ depending only on
$\delta_{0}$ and $\bar{r}$ such that \[
\mathbb{P}\left(\left|f_{n}(u_{0})\right|>\eta_{0}\right)\leq K\exp\left(-\frac{n\beta\eta_{0}^{2}}{c_{1}'\alpha^{2}+2c_{2}'\sqrt{D}\alpha\eta_{0}}\right).\]
Let us set $x_{0}$ such that $\eta_{0}=\alpha\left(\sqrt{c_{1}'\left(x_{0}/\beta\right)}+c'_{2}\sqrt{D}\left(x_{0}/\beta\right)\right)$.
Then: \[
x_{0}\leq\frac{\beta\eta_{0}^{2}}{c_{1}'\alpha^{2}+2c_{2}'\sqrt{D}\alpha\eta_{0}}\]
and \[
\mathbb{P}\left(f_{n}(u_{0})>\eta_{0}\right)\leq K\exp\left(-nx_{0}\right).\]
Then \begin{equation}
P_{1}\leq K\sum_{u_{0}\in\rond{E}_{0}}\exp\left(-nx_{0}\right)\leq K\exp\left(H_{0}-nx_{0}\right).\label{bruitsigma_eq:majoration_P1}\end{equation}
We have that \[
\left\Vert u_{k}-u_{k-1}\right\Vert _{\pi}^{2}\leq2\left(\left\Vert u-u_{k-1}\right\Vert _{\pi}^{2}+\left\Vert u-u_{k}\right\Vert _{\pi}^{2}\right)\leq5\delta_{k-1}^{2}/2\]
then $\left\Vert u_{k}-u_{k-1}\right\Vert _{n}^{2}\leq15\delta_{k-1}^{2}/4$.
As $u_{k-1},u_{k}\in\rond{E}_{k-1}\times\rond{E}_{k}$, it follows
that $\left\Vert u_{k}-u_{k-1}\right\Vert _{\infty}^{2}\leq5\delta_{k-1}^{2}\bar{r}^{2}/2$.
There exists two constants $c_{3}$ and $c_{4}$ such that: \begin{eqnarray*}
\mathbb{P}_{n}\left(\left|f_{n}(u_{k}-u_{k-1})\right|>\eta_{k}\right) & = & \mathbb{P}_{n}\left(\left|f_{n}(u_{k}-u_{k-1})\right|>\eta_{k},\;\left\Vert u_{k}-u_{k-1}\right\Vert _{n}^{2}\leq15\delta_{k-1}^{2}/4\right)\\
 & \leq & K\exp\left(-\frac{n\beta\eta_{k}^{2}}{c_{3}\alpha^{2}\delta_{k-1}^{2}+2c_{4}\alpha\delta_{k-1}}\right).\end{eqnarray*}
Let us fix $x_{k}$ such that $\eta_{k}=\delta_{k-1}a\left(\sqrt{c_{3}\left(x_{k}/\beta\right)}+c_{4}\left(x_{k}/\beta\right)\right)$.
We obtain: \[
x_{k}\leq\frac{\beta\eta_{k}^{2}}{c_{3}\alpha^{2}\delta_{k-1}^{2}+2c_{4}\alpha\delta_{k-1}}\]
and \[
\mathbb{P}\left(\left|f_{n}(u_{k}-u_{k-1})\right|>\eta_{k}\right)\leq K\exp\left(-nx_{k}\right).\]
Then, $P_{2,k}\leq K\exp\left(H_{k-1}+H_{k}-nx_{k}\right)$ and \begin{equation}
P_{2}=\sum_{k=1}^{\infty}P_{2,k}\leq K\sum_{k=1}^{\infty}\exp\left(H_{k-1}+H_{k}-nx_{k}\right).\label{bruitsigma_eq:majoration_P2}\end{equation}
Let us set $\tau>0$ and choose $(x_{k})$ (and then $(\eta_{k})$)
such that \[
\begin{cases}
\sqrt{D}nx_{0}=H_{0}+D+\tau\\
nx_{k}=H_{k-1}+H_{k}+(k+1)D+\tau.\end{cases}\]
Collecting the results, we obtain, by \eqref{bruitsigma_eq:majorationP},
\eqref{bruitsigma_eq:majoration_P1} and \eqref{bruitsigma_eq:majoration_P2}:
\begin{equation}
\mathbb{P}\left(\sup_{u\in\rond{B}_{m,m'}}\left|f_{n}(u)\right|>\eta\right)\leq C\left(e^{-D}e^{-\tau}+e^{-\sqrt{D}}e^{-\tau/\sqrt{D}}\right).\label{bruitsigma_eq:majorationP_D}\end{equation}
It remains to compute $\eta^{2}$. We denote by $C$ a constant depending
only on $\delta_{0}$ and $\bar{r}$ . This constant may vary from
one line to another. We have that:

\begin{eqnarray*}
\eta=\sum_{k=0}^{\infty}\eta_{k} & \leq & C\alpha\left(\sum_{k=1}^{\infty}\delta_{k-1}\left(\sqrt{\frac{x_{k}}{\beta}}+\frac{x_{k}}{\beta}\right)\right)+\alpha\left(\sqrt{\frac{x_{0}}{\beta}}+\sqrt{D}\frac{x_{0}}{\beta}\right).\end{eqnarray*}
Let us recall that $H_{k}=C(k+1)D$. Then, $nx_{k}=C(3k+2)D+\tau$
, $\sqrt{D}nx_{0}=CD+\tau$ and \[
\sum_{k=0}^{\infty}\frac{\delta_{k-1}x_{k}}{\beta}\leq\frac{1}{n\beta}\sum_{k=0}^{\infty}2^{-(k-1)}(C(3k+2)D+\tau)\leq C\frac{D+\tau}{n\beta}.\]
Moreover, \[
\sum_{k=0}^{\infty}\delta_{k-1}\sqrt{\frac{x_{k}}{\beta}}\leq C\frac{\sqrt{D}+\sqrt{\tau}}{\sqrt{n\beta}}.\]
As $D/n\beta\leq1$, there exists a constant $\kappa$ such that \[
\eta^{2}\leq\kappa\alpha^{2}\left(\frac{D}{n\beta}+2\frac{\tau}{n\beta}+\frac{\tau^{2}}{n^{2}\beta^{2}}\right).\]
Then, according to \eqref{bruitsigma_eq:majorationP_D}: \begin{equation}
\mathbb{P}\left(\sup_{u\in\rond{B}_{m,m'}}f_{n}^{2}(u)>\kappa\alpha^{2}\left(\frac{D}{n\beta}+2\frac{\tau}{n\beta}+\frac{\tau^{2}}{n^{2}\beta^{2}}\right)\right)\leq C\left(e^{-D-\tau}+e^{-\sqrt{D}-\tau/\sqrt{D}}\right).\label{sigmabruit_eq:majorationPn}\end{equation}
Furthermore \begin{eqnarray*}
E & := & \E\left(\left[\sup_{u\in\rond{B}_{m,m'}}f_{n}^{2}(u)-\kappa a^{2}\frac{D}{n\beta}\right]_{+}\right)\\
 & = & \int_{0}^{\infty}\mathbb{P}\left(\sup_{u\in\rond{B}_{m,m'}}f_{n}^{2}(u)>\kappa a^{2}\frac{D}{n\beta}+\tau\right)d\tau\end{eqnarray*}
Setting $\tau=\kappa\alpha^{2}\left(2y/n\beta+y^{2}/n^{2}\beta^{2}\right)$,
it follows: \[
E=C\gamma^{2}\int_{0}^{\infty}\mathbb{P}\left(\sup_{u\in\rond{B}_{m,m'}}f_{n}^{2}(u)>\kappa\alpha^{2}\left(\frac{D}{n\beta}+2\frac{y}{n\beta}+\frac{y^{2}}{n^{2}\beta^{2}}\right)\right)\left(\frac{2}{n\beta}+\frac{2y}{n^{2}\beta^{2}}\right)dy.\]
By \eqref{sigmabruit_eq:majorationPn}, \begin{eqnarray*}
E & = & C\kappa\alpha^{2}\left(e^{-D}+e^{-\sqrt{D}}\right)\left(\frac{1}{n\beta}\int_{0}^{\infty}ye^{-y/\sqrt{D}}dy\right)\\
 & \leq & C\frac{\kappa\alpha^{2}}{n\beta}D^{3/2}e^{-\sqrt{D}}.\end{eqnarray*}

\paragraph*{Acknowledgement:}

the author wishes to thank M. Reiss and V. Genon-Catalot for helpful
discussions.

\bibliographystyle{style}
\bibliography{biblio2}

\end{document}